\documentclass[a4paper,11pt,twoside]{article}
\usepackage{pst-fill,pst-grad}
\usepackage{textcomp}  
\usepackage[english]{babel}  
\usepackage[utf8]{inputenc}  
\usepackage[titletoc]{appendix}
\usepackage{titlesec} 
\usepackage{graphicx}  
\usepackage{amsmath} 
\usepackage{float}    
\usepackage{fancyhdr} 
\usepackage[matrix,arrow,curve]{xy} 
\usepackage{pstricks} 
\usepackage{amsmath,amsfonts,verbatim,afterpage,theorem,euscript,mathrsfs,amssymb}
\usepackage{amsfonts}
\usepackage{amssymb}
\usepackage{array} 
\usepackage{dsfont}
\usepackage[colorlinks=true,linkcolor=blue,citecolor=red]{hyperref}
\usepackage{authblk}
\usepackage{color}
\newtheorem{Definition}{Definition}[section] 
\newtheorem{Proposition}{Proposition}[section]

\newtheorem{Lemme}{Lemma}[section]
\newtheorem{Theoreme}{Theorem}[section]

\newtheorem{Corollaire}{Corollary}[section]
\newtheorem{Remarque}{Remark}

\def \vf{\vec{f}} 

\def \vu{\vec{u}}
\def \vv{\vec{v}}
\def \P{\mathbb{P}}

\def \R{\mathbb{R}}

\def \Rt{\mathbb{R}^3}

\def \finpv{\hfill $\blacksquare$} 
\def \pv{{\bf{Proof.}}~} 

\def \ds{\displaystyle}

\title{\bf On the existence, regularity and uniqueness of $L^p$-solutions to the  steady-state 3D Boussinesq system in the whole space}

\author[1]{ Oscar Jarr\'in\footnote{corresponding author: oscar.jarrin@udla.edu.ec}}

\affil[1]{\scriptsize Escuela de Ciencias Físicas y Matemáticas, Universidad de Las Américas, Vía a Nayón, C.P.170124, Quito, Ecuador.} 
\date{\today}
\begin{document}
	\maketitle	 
	\begin{abstract}
 We consider the steady-state Boussinesq system in the whole three-dimensional space,  with the action of external forces and the gravitational acceleration. First, for $3<p\leq +\infty$ we prove the existence of  weak $L^p$-solutions. Moreover, within the framework of a slightly  modified system, we discuss the possibly non-existence of $L^p-$solutions for $1\leq p \leq 3$. Then,  we use the  more general  setting of the $L^{p,\infty}-$spaces to  show that weak solutions and their derivatives are   H\"older continuous functions, where   the maximum gain of regularity is determined by the initial regularity of the external forces and the gravitational acceleration. As a bi-product, we  get a new regularity criterion for the steady-state Navier-Stokes equations. Furthermore, in the particular homogeneous case  when the external forces are equal to zero; and for a range of values of the parameter $p$, we show that weak solutions are not only smooth enough, but also  they are identical to the trivial (zero) solution. This  result is of independent interest, and it is  also known as the Liouville-type problem for the  steady-state Boussinesq system. \\[3mm]
\textbf{Keywords:} Boussinesq system; steady-state $L^p$-solutions; non-existence of solutions; regularity criterion, Lorentz spaces, Liouville problem.   \\[3mm]
\textbf{AMS Classification:} 35A01, 35B53, 35B65.  
	\end{abstract}

\section{Introduction} 
In this note, we consider the steady-state (time-independent) incompressible 3D Boussinesq system in the whole three-dimensional  space. For $\vu : \Rt \to \Rt$ a divergence-free velocity field, $P:\Rt \to \R$ a pressure and $\theta: \Rt \to \R$ the temperature of the fluid we have:
\begin{equation}\label{Boussinesq}
\begin{cases}\vspace{2mm}
-\Delta \vu + \text{div}(\vu \otimes \vu) + \vec{\nabla} P = \theta \vec{{\bf g}}+ \vf,  \qquad \text{div}(\vu)=0, \\
-\Delta \theta + \text{div}(\theta \vu) =  g,
\end{cases}  \qquad \text{in} \quad \Rt, 
\end{equation} 
where $\vec{{\bf g}}: \Rt \to \Rt$ is the gravitational acceleration  and $\vf : \Rt \to \Rt$, $g : \Rt \to \R$ denote external forces acting on each equation of this system. Moreover, with a minor loss of generality all the physical constants have been set to be  equal to one.  

\medskip

The system (\ref{Boussinesq})  describes the dynamics of a viscous incompressible fluid with heat exchanges \cite{Boussinesq,Chandrasekhar}, and moreover, it takes into account 
the physical phenomena carried out by a gravitational field $\vec{\bf g}$, like the Earth's gravitational field.   This system arises from an approximation on a system coupling the classical Navier-Stokes equations and the equations of thermodynamics. In this approximation, the variations of the density due to heat transfers are neglected in the continuity equation, but are taken into account in the equation of the motion through an
additional buoyancy term proportional to the temperature variations and the gravitational acceleration: $\theta \vec{{\bf g }}$. In our study, we shall observe that this term yields interesting qualitative properties to the system (\ref{Boussinesq}), in particular the possible non-existence of $L^p$-solutions for some range of values of the parameter $p$. 

\medskip

When we set $\theta \equiv 0$, the system (\ref{Boussinesq}) boils down to the classical steady-state Navier-Stokes equations: 
\begin{equation}\label{NS}
-\Delta \vu + \text{div}(\vu \otimes \vu) + \vec{\nabla} P =  \vf, \qquad \text{div}(\vu)=0, \qquad \text{in} \ \ \Rt,
\end{equation}
which contains many challenging open problems. We refer to the books \cite{PGLR}  and \cite{PGLR1} for a complete  theoretical study of these equations.  

\medskip

In the \emph{homogeneous case} when $\vf \equiv 0$ and $g\equiv 0$, the steady state Boussinesq   system  has been mainly  studied in the setting of a \emph{bounded and  smooth  domain} $\Omega \subset \Rt$, jointly with the Dirichlet conditions:
\begin{equation}\label{Boussinesq-Domain}
\begin{cases}\vspace{2mm}
- \nu\Delta \vu + \text{div}(\vu \otimes \vu) + \vec{\nabla} P = \theta \vec{{\bf g}}, \qquad \text{div}(\vu)=0, \quad \text{in} \ \ \Omega, \\ \vspace{2mm} 
-\Delta \theta + \text{div}(\theta \vu) = 0, \quad \text{in} \ \ \Omega, \\
\vu = \vu_b, \quad \theta=\theta_b \quad \mbox{on} \ \ \partial \Omega.
\end{cases}
\end{equation}
Here, $\nu>0$ is the viscosity parameter, and moreover, the given functions $\vu_b: \partial \Omega \to \Rt$ and $\theta_b : \partial \Omega \to \R$ are  boundary data.  The qualitative properties of this system: existence, uniqueness and regularity of solutions, deeply depend on suitable hypothesis on the functions $\vec{\bf g}$, $\vu_b$ and $\theta_b$ as well as on suitable hypothesis on  the boundary $\partial \Omega$. 

\medskip 

When  $\partial \Omega$ is assumed to be a  \emph{connected set} of class $\mathcal{C}^2$, and additionally, when the viscosity parameter $\nu$ is large enough,  in \cite{Santos,Villamizar} it is proven the existence and uniqueness  of weak solutions $(\vu, \theta) \in L^3(\Omega) \times L^2(\Omega)$ to the system (\ref{Boussinesq-Domain}) with boundary data $(\vu_b, \theta_b) \in L^2(\partial \Omega) \times L^2(\partial \Omega)$ and $\vec{{\bf g}} \in L^{\infty}(\Omega)$.  

\medskip

 Thereafter, the case of more irregular boundary data was studied in \cite{Kim}. For  $\vec{{\bf g}} \in L^{\infty}(\Omega)$, the existence and uniqueness of weak solutions $(\vu, \theta) \in L^{p}(\Omega) \times L^q(\Omega)$ to the system (\ref{Boussinesq-Domain}) are proven under the  smallness  assumption $ \| \vu_b \|_{W^{-1/p,p}(\partial \Omega)}  + \|  \theta_b \|_{W^{-1/q,q}(\partial \Omega)} \ll 1$; and  where the parameters $p$ and $q$ verify the following set of technical conditions: 
 with $3 \leq p <+\infty$,  $ p' < q < +\infty$ (where $1/p +1/p'=1$) and   $3p / (3 + 2p) \leq r$. 

\medskip

On the other hand, the more general case when $\partial \Omega$ is of class $\mathcal{C}^{1,1}$ but \emph{not necessarily connected} is addressed in \cite{Avecedo}. In this work, for a weaker assumption on  gravitational acceleration: $\vec{{\bf g}} \in L^{3/2}(\Omega)$,   it is proven the existence of weak solutions in the space $H^1(\Omega)$ just considering smallness of $\vu_b$ across each connected component $\partial \Omega _i$ of the boundary $\partial \Omega$. Moreover,  it is also proven a gain of regularity of weak solutions in the space $W^{1,p}(\Omega)$ with $p>2$; and in the space $W^{2,p}(\Omega)$ with $p\geq 6/5$.  For more related results, we refer to \cite{Gil,Maimoto1,Marimoto2} and the references therein. 

\medskip 

To the best of our knowledge, these qualitative properties of $L^p$-solutions to the system (\ref{Boussinesq}) (in  whole space $\Rt$) have not been studied before. It is worth highlighting  the previously cited works are not longer valid in $\Rt$ due to the lack of some of their key tools,  for instance, the  embedding properties of the $L^p(\Omega)$- spaces and the \emph{compact} Sobolev embeddings.  We thus use different approaches to study the existence, regularity and uniqueness of $L^p$-solutions to the system  (\ref{Boussinesq}); and our key tools are mainly based in useful properties of the Lorentz spaces.  We refer to the book \cite{DCh} for a complete study of these spaces.

\medskip

We recall that for a measurable function  $f: \Rt \to \R$  and for a parameter $\lambda \geq 0$ we define  the distribution function
\[ d_f(\lambda)= dx \left( \left\{  x \in \Rt: \ | f(x)|>\lambda \right\} \right), \]
where $dx$ denotes the Lebesgue measure. Then, the re-arrangement function $f^{*}$ is defined by the expression 
\[ f^{*}(t)= \inf \{ \lambda \geq 0: \ d_f(\lambda) \leq t\}. \]
By definition, for $1\leq p <+\infty$ and $1\leq q \leq +\infty$ the Lorentz space $L^{p,q}(\Rt)$ is the space  of measure functions $f: \Rt \to \R$ such that $\| f \|_{L^{p,q}}<+\infty$, where:
\begin{equation*}
\| f \|_{L^{p,q}}=
\begin{cases}\vspace{3mm}
\ds{\frac{q}{p} \left( \int_{0}^{+\infty} (t^{1/p} f^{*}(t))^{q} dt\right)^{1/q}}, \quad q<+\infty, \\
\ds{\sup_{t>0} \, t^{1/p} f^{*}(t)}, \quad q=+\infty.
\end{cases}
\end{equation*}
It is worth mentioning some important properties of these spaces. The quantity $\| f \|_{L^{p,q}}$ is often used as a norm, even thought it does not verify the triangle inequality. However, there exists  an equivalent norm (strictly speaking) which makes these spaces into Banach spaces. On the other hand, these spaces are homogeneous of degree $-3/p$ and for $1 \leq q_1 <  p < q_2\leq +\infty$ we have the continuous embedding 
\[ L^{p,q_1}(\Rt) \subset L^{p}(\Rt)=L^{p,p}(\Rt) \subset L^{p,q_2}(\Rt). \]  
Finally, for $p=+\infty$ we also have the identity  $L^{\infty,\infty}(\Rt)=L^\infty(\Rt)$.

\medskip 

In our first main result, we find some  smallness conditions on the external forces and the gravitational acceleration to construct $L^p$-weak solutions to the system (\ref{Boussinesq}). Precisely, by following some of the ideas in \cite{Bjorland}, first we solve this system in the $L^{3,\infty}-$space, where (to the best of our knowledge) we are able to apply a fixed point argument. Thereafter,  we show that these $L^{3,\infty}-$solutions also belong to the space $L^p(\Rt)$ (with $3< p  \leq +\infty$) as long as the external forces  verify additionally suitable conditions. 

\begin{Theoreme}\label{Th1}   Assume that the external forces $\vf$ and $g$ satisfy  
\[ (-\Delta)^{-1}\P (\vf) \in L^{3,\infty}(\Rt) \quad \mbox{and} \quad  (-\Delta)^{-1} g \in L^{3,\infty}(\Rt). \]	
Moreover, assume that the gravitational acceleration $\vec{{\bf g}}$ satisfies
\[ \vec{{\bf g}} \in L^{3/2,1}(\Rt). \]
There exists a universal quantity $\varepsilon>0$ such that if 
\begin{equation*}
\delta=\left\| (-\Delta)^{-1}\P (\vf)   \right\|_{L^{3,\infty}} + \left\| (-\Delta)^{-1} g \right\|_{L^{3,\infty}}< \varepsilon \quad \mbox{and} \quad \| \vec{\bf g} \|_{L^{3/2,1}} < \varepsilon,
\end{equation*}	
then the following statements hold: 
\begin{enumerate}
	\item The system (\ref{Boussinesq}) has a weak solution $(\vu, \theta) \in L^{3,\infty}(\Rt)$ satisfying and uniquely defined by the condition $\| \vu \|_{L^{3,\infty}}+\| \theta \|_{L^{3,\infty}} \leq 3 \delta$.
		
	\item This solution verifies $(\vu,\theta) \in L^p(\Rt)$, with $3<p\leq  +\infty$, as long as 
	\begin{equation*}
	\begin{cases}\vspace{2mm}
	(-\Delta)^{-1}  \P(\vf)  \in L^{p} \cap L^4 (\Rt), \ \  (-\Delta)^{-1}g \in L^{p} \cap L^4 (\Rt), \ \  \mbox{when}\ \ 3<p<4, \\
	 (-\Delta)^{-1}\P(\vf)  \in L^{p}  (\Rt), \ \ (-\Delta)^{-1}g \in L^{p} (\Rt), \ \ \mbox{when} \ \  4 \leq p \leq +\infty,
	\end{cases} \ \ \mbox{and} \ \ \ \vec{\bf g} \in  L^{3,1}(\Rt).
	\end{equation*}
\end{enumerate}	
\end{Theoreme}	

Some comments are in order. As mentioned in the first point above,   the smallness assumption $\delta <\varepsilon$  yields the existence of weak $L^{3,\infty}$-solutions, where this particular space naturally appears in the involved fixed point estimates (see the Proposition \ref{Prop-Existence-L3} below for more details). However, the  conditions $(-\Delta)^{-1}\P (\vf) \in L^{3,\infty}(\Rt)$ and $ (-\Delta)^{-1} g \in L^{3,\infty}(\Rt)$ allow us to consider very rough external forces, for instance, for suitable constants $\eta_1,\eta_2>0$ and for the Dirac mass $\delta_0$,  we can consider $\vf = \eta_1 (\delta_0, \delta_0, \delta_0)$ and $g = \eta_2 \delta_0$.   On the other hand, the gravitational acceleration $\vec{\bf g}$ also plays an important role in this existence result. In particular, the required condition $\vec{\bf g} \in L^{3/2,1}(\Rt)$ is similar to the  condition  $\vec{\bf g} \in L^{3/2}(\Omega)$ found in \cite{Avecedo}. 

\medskip

The second point above shows a persistence property of solutions in the $L^p$-spaces (see the Proposition \ref{Prop-Existence-Lp} below for all the details). Under minor technical modifications,  this result is also valid in the (more general) setting of the $L^{p,q}-$spaces.  

\medskip

As may observe, this last  result is valid for the range of values  $3<p\leq +\infty$ and it is natural to look  for similar   persistence properties in the complementary range of values $1\leq p \leq 3$. It is thus interesting to remark that the constraint $3<p\leq +\infty$ is  not only  technical  (see the Remarks \ref{Rmk-2} and \ref{Rmk3} below),  but it seems to be a phenomenological effect of the system (\ref{Boussinesq}). More precisely, the term $\theta \vec{\bf g}$ in the first equation of this system potentially  yields the non-existence of $L^p$-solutions for $1\leq p \leq 3$. 

\medskip

To study this  phenomena, let us consider the following toy model for the system (\ref{Boussinesq}): 
\begin{equation}\label{Boussinesq-toy-model}
\begin{cases}\vspace{2mm}
-\Delta \vu + \text{div}(\varphi \vu \otimes \vu) =  \theta \vec{{\bf g}} + \vf, \quad \text{div}(\vu=0), \qquad |\varphi(x)| \leq \frac{C}{1+|x|^2}, \\
-\Delta \theta + \text{div}(\theta \vu)= g. 
\end{cases}
\end{equation} 
Here, the localizing function $\varphi$ in the nonlinear term is a merely technical requirement  and, in further research, we aim to get rid of this function. However, the main goal  of this toy model is to highlight  the effects of the term $\theta \vec{{\bf g}}$ in the non-existence of solutions. 

\medskip

Since  $\varphi \in L^{\infty}(\Rt)$ the  modified bilinear  term $\text{div}(\varphi \vu \otimes \vu)$   follows the same estimates of the original one in the equation (\ref{Boussinesq}). See the Section \ref{Sec:Proof-Th1} for all the details.  Consequently,  the result obtained in Theorem \ref{Th1} is valid for the system (\ref{Boussinesq-toy-model}). In particular, one can consider external forces $\vf, g \in \mathcal{S}(\Rt)$ and a gravitational acceleration $\vec{{\bf g}} \in \mathcal{S}(\Rt)$   verifying 
\begin{equation}\label{Force-Supp-Fourier}
0 \notin \text{supp}(\widehat{\vf}), \quad  0 \notin \text{supp}(\widehat{g}) \quad \mbox{and} \quad 0 \notin \text{supp}(\widehat{\vec{\bf {g}}}),
\end{equation} 
hence we have 
\begin{equation*}
\ds{(-\Delta)^{-1} \vf\in \mathcal{S}(\Rt)}, \quad  \ds{(-\Delta)^{-1} g \in \mathcal{S}(\Rt)} \quad \mbox{and} \quad  (-\Delta)^{-1} \vec{{\bf g}} \in \mathcal{S}(\Rt).
\end{equation*}
Then, we assume the smallness condition in  Theorem \ref{Th1} and  there exists a solution $(\vu, \theta) \in L^{3,\infty}(\Rt)$ to the system (\ref{Boussinesq-toy-model}) which verifies $(\vu, \theta) \in L^p(\Rt)$ with $3<p\leq +\infty$. In this context, our next result reads as follows:
\begin{Proposition}\label{Prop-Nonexistence} 
There exist    $\vf, g, \vec{{\bf g}} \in \mathcal{S}(\Rt)$ well-prepared functions satisfying (\ref{Force-Supp-Fourier}), such that  the associated solution $(\vu,\theta)$ to the system  (\ref{Boussinesq-toy-model}) verifies $ \vu \notin L^p(\Rt)$ for $1 \leq p \leq 3$.
\end{Proposition}	

The proof  is essentially  based on the following estimate from below (see the Proposition \ref{Prop-Estimate-Below} for more details)
\[  \frac{1}{|x|}  \left| \int_{\Rt} \theta(y)\vec{\bf g}(y) dy \right| \lesssim | \vu(x)|, \qquad |x| \gg 1, \]
which yields $ \vu \notin L^p(\Rt)$ for the range of values  $1 \leq p \leq 3$. Here, the well-prepared data $\vf, g$ and $\vec{\bf g}$ ensure that $\ds{\int_{\Rt} \theta(y)\vec{\bf g}(y) dy \neq 0}$.  It is also interesting to compare this result to the one obtained in \cite[Theorem $3.5$]{Bjorland} for the steady-state Navier-Stokes equations (\ref{NS}). In this work, it is proven a non-existence  of  $L^p$-solutions for $1\leq p \leq 3/2$ and we thus observe the strongest effects of the term $\theta \vec{\bf g}$ is the Boussinesq system (\ref{Boussinesq}). 

\medskip

In our second main result, we study another persistence problem associated with the system (\ref{Boussinesq}), which  is an important question to a better mathematically comprehension of this system:  the gain of regularity of weak solutions from an initial regularity in the data $\vf, g$ and $\vec{\bf g}$.  We shall consider here a fairly general notion of weak solutions, which is given in the following: 
\begin{Definition}\label{Def-weakSol-Boussinesq} Let $\vf, g \in \mathcal{D}'(\Rt)$ and let $\vec{{\bf g}} \in L^{2}_{loc}(\Rt)$. A  weak solution of the coupled system (\ref{Boussinesq}) is  a triplet $(\vu, P, \theta)$,  where:  $\vu \in L^{2}_{loc}(\Rt)$, $P \in \mathcal{D}'(\Rt)$  and $\theta \in L^{2}_{loc}(\Rt)$, such that  it verifies (\ref{Boussinesq})  in the distributional sense. 
\end{Definition}	
It is worth observing we use minimal conditions on the functions $\vf,g,\vec{\bf g},\vu, P$ and $\theta$  to ensure that  all the terms in  the system  (\ref{Boussinesq})  are well defined as distributions. In addition, we let the pressure $P$ to be a very general object since  we only assume $P\in \mathcal{D}'(\Rt)$. 

\medskip

The $L^{p,\infty}$-spaces provide us a general and suitable framework to improve the regularity of weak solutions defined above.  For  $1\leq p <+\infty$  and for a regularity parameter $k\in \mathbb{N}$, we introduce the Sobolev-Lorentz space 
\begin{equation*}
\mathcal{W}^{k,p}(\Rt)=\left\{  f \in L^{p,\infty}(\Rt): \ \partial^{\alpha} f \in L^{p,\infty}(\Rt) \  \mbox{for all multi-indice} \ |\alpha| \leq k  \right\}. 
\end{equation*}
Moreover,  we denote by $W^{k,\infty}(\Rt)$ the classical Sobolev space of bounded  functions with bounded  weak derivatives until the order $k$. Finally, for  $0<s<1$ we shall denote by $\mathscr{C}^{k,s}(\Rt)$ the H\"older space of $\mathcal{C}^{k}-$ functions  whose   derivatives are H\"older continuous functions with parameter $s$. In this setting, our second result reads as follows:

\begin{Theoreme}\label{Th2} Let $(\vu, P, \theta)$ be a weak solution to the coupled system (\ref{Boussinesq}) given by Definition \ref{Def-weakSol-Boussinesq}. Assume that 
\begin{equation*}
\vu \in L^{p,\infty}(\Rt) \quad \mbox{and} \quad \theta \in  L^{p,\infty}(\Rt), \quad \mbox{with} \ \ 3<p<+\infty. 
\end{equation*}	
Then, if for $k \geq 0$: 
\begin{equation}\label{Condition-Reg-force}
\vf , \, g  \in \mathcal{W}^{k,\frac{3p}{3+p}} \cap W^{k,\infty}(\Rt),  
\end{equation}	 
and 
\begin{equation}\label{Cond-Reg-g}
\vec{\bf {g}}  \in  \mathcal{W}^{k,\frac{3p}{3+p}} \cap W^{k,\infty}(\Rt),
\end{equation}
it follows that $\vu \in \mathcal{W}^{k+2,p}(\Rt)$, $P\in  \mathcal{W}^{k+1,p}(\Rt)$ and $\theta \in \mathcal{W}^{k+2,p}(\Rt)$. Moreover, we have  
\[  \vu \in \mathscr{C}^{k+1,s}(\Rt), \ \  P\in \mathscr{C}^{k,s}(\Rt)  \ \ \mbox{and} \ \  \theta \in  \mathscr{C}^{k+1,s}(\Rt), \ \ \mbox{with} \ \  s=1-3/p. \]
\end{Theoreme}

In (\ref{Condition-Reg-force}) and (\ref{Cond-Reg-g})  the parameter $k$ measures the  initial regularity of the data, which yields a gain of regularity of weak solutions of the order $k+2$. This (expected) maximum gain of regularity is given by the effects of the Laplacian operator in both equations of the system (\ref{Boussinesq}). We refer to Remark \ref{Rmk4} for more details. 

\medskip

Let us briefly explain the general strategy of the proof, which   bases on two key ideas. First, by assuming $\vu \in L^{p,\infty}(\Rt)$ and $\theta \in L^{p,\infty}(\Rt)$ and by using useful properties of the parabolic (time-dependent) Boussinesq system (
given in the equation (\ref{Boussinesq-Parabolic})), we prove that $\vu$ and $\theta$ are bounded functions on $\Rt$. Thereafter, we use a bootstrap argument to show that    $\vu \in \mathcal{W}^{k+2,p}(\Rt)$  and $\theta \in \mathcal{W}^{k+2,p}(\Rt)$. Finally, we use  the more general setting of the Morrey spaces to conclude the wished H\"older regularity.  

\medskip

These ideas can also be applied to the other relevant coupled system of the fluid dynamics, for instance, the \emph{magneto-hydrodynamics equations} and some variations of the Boussinesq system as the \emph{Bérnard system} and the \emph{magnetic-Bérnard system}. In particular,  when we set $\theta \equiv 0$ the system (\ref{Boussinesq}) agrees with the Navier-Stokes equations (\ref{NS}); and as a  direct consequence of the theorem above  we obtain a new regularity criterion for these equations: 
\begin{Corollaire} Let $(\vu, P) \in L^{2}_{loc}(\Rt) \times \mathcal{D}'(\Rt)$ be a weak solution to the equation  (\ref{NS}). If $\vu \in L^{p,\infty}(\Rt)$ with $p>3$ and if the external force $\vf$ verifies (\ref{Condition-Reg-force}) for some $k\geq 0$, then we have $\vu \in \mathcal{W}^{k+2,p}(\Rt)$ and $P\in \mathcal{W}^{k+1,p}(\Rt)$. In addition, we have $\vu \in \mathscr{C}^{k+1,s}(\Rt)$ and $P\in \mathscr{C}^{k,s}(\Rt)$ with $s=1-3/p$.
\end{Corollaire}

Another interesting  consequence of Theorem \ref{Th2} arises when we consider  the homogeneous  system (when $\vf\equiv 0$ and $g\equiv 0$):
\begin{equation}\label{Boussinesq-homog}
\begin{cases}\vspace{2mm}
-\Delta \vu + \text{div}(\vu \otimes \vu) + \vec{\nabla} P = \theta \vec{{\bf g}}, \qquad \text{div}(\vu)=0, \\
-\Delta \theta + \text{div}(\theta \vu) = 0,
\end{cases}  \qquad \text{in} \quad \Rt,
\end{equation}
and then, the gain of regularity of weak solution is now determined by the gravitational acceleration  $\vec{\bf g}$ as long as it verifies (\ref{Cond-Reg-g}). In particular, we have $(\vu, P, \theta) \in \mathcal{C}^{\infty}(\Rt)$, provided that $\vu, \theta \in L^{p,\infty}(\Rt)$ with $3<p<+\infty$ and $\ds{\vec{\bf g} \in \bigcap_{k \geq 0} \mathcal{W}^{k,\frac{3p}{3+p}} \cap W^{k,\infty}(\Rt)}$.

\medskip

Regularity of weak solutions to the homogeneous system (\ref{Boussinesq-homog}) is also one of the key assumptions when studying another relevant problem: the uniqueness of the trivial solution $(\vu, P, \theta)=(0,0,0)$. This is problem is commonly known as the \emph{Liouville-type  problem} for the homogeneous  Boussinesq system (\ref{Boussinesq-homog}) and, to the best of our knowledge, it has not been studied before. Thus, in our third main result, we find a range of values of the parameter $p$ for which  $L^{p,q}(\Rt)$-solutions to the system (\ref{Boussinesq-homog}) vanish identically.

\begin{Theoreme}\label{Th3}  Let $(\vu, P, \theta)$ be a weak solution to the coupled homogeneous  system (\ref{Boussinesq-homog}) given by Definition \ref{Def-weakSol-Boussinesq}. The trivial solution $(\vu, P, \theta)=(0,0,0)$ in uniquely determined in the following cases:
\begin{enumerate}
\item When $\vu \in L^{\frac{9}{2}, q}(\Rt)$, $\theta \in L^{\frac{9}{2},q}(\Rt)$ with $1\leq q<+\infty$, and $\vec{\bf g} \in \mathcal{W}^{1,9/2} \cap W^{1,\infty}(\Rt)$. 
\item When $\vu \in L^{p,\infty}(\Rt)$, $\theta \in  L^{p,\infty}(\Rt)$ and $\vec{\bf g} \in \mathcal{W}^{1,p} \cap W^{1,\infty}(\Rt)$, with $3<p < \frac{9}{2}$. 
\end{enumerate}		 
\end{Theoreme}

The proof of this theorem is  based on fine local estimates for each equation in the system (\ref{Boussinesq-homog}), which are known as the \emph{Caccioppoli-type estimates} (see the Proposition \ref{Caccioppoli} below). To do this, a minimal regularity of solutions is required,  more precisely,  we need  that $\vu, \theta \in C^{2}_{loc}(\Rt)$ and $P \in \mathcal{C}^{1}_{loc}(\Rt)$. Thus, this regularity is ensured by the assumptions on $\vu, \theta$ and $\vec{\bf g}$ given above; and by a direct application of Theorem \ref{Th2}. 
 
\medskip

The Liouville-type problem has been extensively studied in the homogeneous Navier-Stokes equations (\ref{NS}) (when $\vf \equiv 0$). See, for instance, \cite{ChaeYoneda,ChaeWolf,ChaeWeng,ChJaLem,Ser2015,Ser2016,SerWang} and the references therein. In this sense, this result is a generalization of these previous works to the coupled framework of the system (\ref{Boussinesq-homog}). In particular, the first point above generalizes the well-known Galdi's result \cite[Theorem $X.9.5$]{Galdi}. Here, it is important to emphasize that the Liouville problem for the equation (\ref{NS}) in the largest space $L^{9/2,\infty}(\Rt)$ and in the setting of the Lebesgue or Lorentz spaces with parameter $p>9/2$ is still an open question far from obvious. The second point above generalizes to the system (\ref{Boussinesq-homog}) some previous results obtained in \cite[Theorem $1$]{Jarrin} for the equation (\ref{NS}). Here, the constraint $3<p<9/2$ allows us to work in the largest space $L^{p,\infty}(\Rt)$. 

\medskip

Finally, let us mention that in future research  some recent results on the Liouville problem for the equation (\ref{NS}), which consider more sophisticated functional spaces, could be also adapted to the system (\ref{Boussinesq-homog}), but this come out the scope of this note. 

\medskip

{\bf Organization of the paper.}   In Section \ref{Sec-Preliminaries} we summarize some useful properties of the $L^{p,q}$-spaces. Section \ref{Sec:Proof-Th1} and \ref{Sec-Non-existence} are devoted to the proof of Theorem \ref{Th1} and Proposition \ref{Prop-Nonexistence} respectively. Then, in Section \ref{Sec-Regularity} we give a proof of Theorem \ref{Th2}; and this note finishes with the proof of Theorem \ref{Th3} given in Section \ref{Sec-Liouville}.

\section{Preliminaries}\label{Sec-Preliminaries} 
For the completeness of this note,   we  summarize  here some  well-known results which will be useful in the sequel.  In the forthcoming estimates, $C>0$ is  a generic constant.   
\begin{Lemme}[Young inequalities]\label{Lem-Young}
Let $1<p,p_1,p_2<+\infty$ and $1\leq q,q_1, q_2 \leq+\infty$. There exists a constant $C_i>0$, which depends on the parameters above, such that  the following estimates hold:  
\begin{enumerate}
	\item $\ds{\Vert f \ast h \Vert_{L^{p,q}} \leq C_1 \Vert f \Vert_{L^{p_1, q_1}}\, \Vert h \Vert_{L^{p_2, q_2}}}$,  with $1+\frac{1}{p}=\frac{1}{p_1}+\frac{1}{p_2}$, $\frac{1}{q}\leq \frac{1}{q_1}+\frac{1}{q_2}$ and   $C_1=C\,p\left( \frac{p_1}{p_1-1}\right)\left( \frac{p_2}{p_2-1}\right)$.  
	\item $\ds{\Vert f \ast h \Vert_{L^{p,q}} \leq C_2\,  \Vert f \Vert_{L^1} \, \Vert h \Vert_{L^{p,q}}}$,  with $C_2= C\, \frac{p^2}{p-1}$.  
	\item $\ds{\Vert f \ast h \Vert_{L^\infty} \leq C_3\,  \Vert f \Vert_{L^{p,q}}\, \Vert h \Vert_{L^{p',q'}}}$, with  $1=\frac{1}{p}+\frac{1}{p'}$, $1 \leq \frac{1}{q}+\frac{1}{q'}$  and $C_3= C \left(\frac{p}{p-1} \right)\left(\frac{p'}{p' -1} \right)$.
\end{enumerate}	 
\end{Lemme}	
For a proof see \cite[Section $1.4.3$ ]{DCh}.
\begin{Lemme}[H\"older inequalities]\label{Lem-Holder} 
\begin{enumerate}
\item[]
\item  	Let $f \in L^{p_1,q_1}(\Rt)$ and $h \in L^{p_2,q_2}(\Rt)$ with $1\leq p_1,p_2<+\infty$, $1\leq q_1,q_2\leq +\infty$ verifying the relationships $1=1/p_1 + 1/p_2$ and $1=1/q_1+1/q_2$. Then we have $fh\in L^1(\Rt)$ and it holds $\ds{\| fh \|_{L^1} \leq \,C \| f \|_{L^{p_1,q_1}}\, \| h \|_{L^{p_2,q_2}}}$.
\item Moreover, for $1\leq p,p_1,p_2<+\infty$, $1\leq q,q_1,q_2 \leq +\infty$ and with the  relationships  $1/p=1/p_1 + 1/p_2$ and $1/q=1/q_1+1/q_2$; we have $fh\in L^{p,q}(\Rt)$ and it holds  $\ds{\| fh \|_{L^{p,q}} \leq C  \| f \|_{L^{p_1,q_1}}\, \| h \|_{L^{p_2,q_2}}}$. 
\end{enumerate}		
\end{Lemme}	
A proof can be found in \cite[Theorems $1.2.6$ and $1.4.1$]{DCh}.

\begin{Lemme}[Interpolation inequalities]\label{Lem-Interpolation} 
\begin{enumerate}
\item[]
\item 		Let $f \in L^{p,\infty}(\Rt)  \cap L^{\infty}(\Rt)$, with $1\leq p < +\infty$.  Then, for all $1\leq \sigma < +\infty$ we have $f \in L^{p\sigma}(\Rt)$  and the following estimate holds:  $\ds{\Vert f \Vert_{L^{p\sigma}}\leq C\, \Vert f \Vert^{\frac{1}{\sigma}}_{L^{p,\infty}}\, \Vert f \Vert^{1-\frac{1}{\sigma}}_{L^{\infty}}}$. 

\item  Let $f \in L^{p_1,q_1}(\Rt) \cap L^{p_2,q_2}(\Rt)$ where $1<p_1<p<p_2<+\infty$ and $1\leq q,q_1,q_2 \leq +\infty$. Then we have $f \in L^{p,q}(\Rt)$ and it holds: $\ds{\| f \|_{L^{p,q}} \leq  C \| f \|^{\sigma_1}_{L^{p_1,q_1}}\, \| f \|^{1-\sigma_1}_{L^{p_2,q_2}}}$, where $0<\sigma_1<1$ depends on $p,p_1$ and $p_2$. 
\end{enumerate}		
\end{Lemme}	
For a proof see \cite[Theorem $1.1.2$ ]{DCh} and \cite[Proposition $2.3$]{PGLR}.  

\begin{Lemme}[Lorentz-Besov embedding]\label{Lem-Lorentz-Besov} For all $t>0$ it  holds:  $\ds{t^{\frac{3}{2p}} \Vert e^{t\Delta} f \Vert_{L^{\infty}} \leq C\, \Vert f \Vert_{L^{p,\infty}}}$, with $1\leq p <+\infty$. 
\end{Lemme} 	
This estimate is a direct consequence of the continuous embedding $L^{p,\infty}(\Rt)\subset \dot{B}^{-\frac{3}{p}, \infty}_{\infty}(\Rt)$. See \cite[Page $171$]{PGLR1}. We recall that the homogeneous Besov space $\dot{B}^{-\frac{3}{p}, \infty}_{\infty}(\Rt)$ can be characterized  as the space of temperate distributions $f\in \mathcal{S}'(\Rt)$ such that $\ds{\sup_{t>0}\, t^{\frac{3}{2 p}}\Vert e^{t\Delta} f \Vert_{L^{\infty}} <+\infty}$.
\begin{Lemme}[Riesz transforms]\label{Lem-Riesz} For $i=1,2,3$ let  $\ds{\mathcal{R}_i = \frac{\partial_i}{\sqrt{-\Delta}}}$ be the i-th Riesz transform. Then, for $i,j=1,2,3$ the operator $\mathcal{R}_{i}\mathcal{R}_{j}$ is continuous in the space $L^{p,\infty}(\Rt)$ (with $1<p<+\infty$) and we have $\ds{\Vert \mathcal{R}_i \mathcal{R}_j (f) \Vert_{L^{p,\infty}} \leq C \Vert f  \Vert_{L^{p,\infty}}}$.  In particular, the Leray's projector $\P$ is a bounded operator in the space $L^{p,\infty}(\Rt)$.
\end{Lemme}	 
A proof can be found in \cite[Lemma $4.2$]{Kato}. Finally, we shall use the following result linking the Morrey spaces and the H\"older regularity of functions.  For a proof see   \cite[Proposition $3.4$]{GigaMiyakawa}. Recall that for $1\leq p < +\infty$ the homogeneous  Morrey  space $\dot{M}^{1,p}(\Rt)$ is defined as the space of locally finite Borel measures $d\mu$ such that 
\begin{equation}\label{Morrey}
\sup_{x_0 \in \Rt,\, R>0} R^{\frac{3}{p}} \left( \frac{1}{\vert B(x_0, R)\vert} \int_{B(x_0,R)} d\vert \mu \vert (x) \right)<+\infty. 
\end{equation}
\begin{Lemme}[H\"older regularity]\label{Lem-Holder-Reg} Let $f \in \mathcal{S}'(\Rt)$ such that $\vec{\nabla} f \in \dot{M}^{1,p}(\Rt)$, with $p>3$. There exists a constant $C>0$ such that for all  $x,y\in \Rt$ we have $\ds{\vert f(x)-f(y) \vert \leq C\, \Vert \vec{\nabla} f \Vert_{\dot{M}^{1,p}} \, \vert x-y \vert^{1-3/p}}$. 
\end{Lemme}

\section{Proof of Theorem \ref{Th1}}\label{Sec:Proof-Th1}
First note that the system (\ref{Boussinesq}) can be rewritten  as the following (equivalent)  problem:
\begin{equation}\label{Boussinesq-P}
\begin{cases}\vspace{2mm}
-\Delta \vu + \P (\text{div}(\vu \otimes \vu))  =  \P(\theta \vec{{\bf g } })+ \P(\vf), \qquad \text{div}(\vu)=0, \\
-\Delta \theta + \text{div}(\theta \vu) = g,
\end{cases} 
\end{equation} 
hence we obtain  the  integral formulations: 
\begin{equation}\label{Boussineq-Int-1}
\ds{\vu  = - \frac{1}{-\Delta} \left( \P ( \text{div} (\vu \otimes \vu) ) \right)+ \frac{1}{-\Delta}\Big( \P \left(  \theta \vec{\bf g} \right)\Big) +  \frac{1}{-\Delta}\Big( \P(\vf )\Big)},  
\end{equation}
\begin{equation}\label{Boussineq-Int-2}
\ds{\theta =- \frac{1}{-\Delta} \left( \text{div}(\theta \vu) \right) +  \frac{1}{-\Delta}( g )}. 
\end{equation}  
In the next propositions, we shall prove each point stated in Theorem \ref{Th1}.
\begin{Proposition}\label{Prop-Existence-L3} Assume that $\ds{\frac{1}{-\Delta}\Big( \P (\vf)\Big) \in L^{3,\infty} (\Rt)}$, $\ds{\frac{1}{-\Delta}(g) \in L^{3,\infty} (\Rt)}$ and $\vec{\bf {g}} \in L^{3/2,1}(\Rt)$. There exists a universal constant $\varepsilon_1>0$ (defined in (\ref{epsilon1})) such that if 
\begin{equation}\label{Control-Forces-Grav}
\delta=\left\|  \frac{1}{-\Delta}\Big( \P(\vf) \Big)\right\|_{L^{3,\infty}} + \left\|  \frac{1}{-\Delta}(g)\right\|_{L^{3,\infty}}< \varepsilon_1 \quad \mbox{and} \quad \| \vec{\bf g} \|_{L^{3/2,1}} < \varepsilon_1,
\end{equation}	
then the system (\ref{Boussineq-Int-1})-(\ref{Boussineq-Int-2}) has a solution $(\vu, \theta) \in L^{3,\infty}(\Rt)$ satisfying and uniquely defined by $\| \vu \|_{L^{3,\infty}}+ \| \theta \|_{L^{3,\infty}} \leq 3 \delta$. 
\end{Proposition}	
\pv To prove this proposition, we shall use the next version of the Picard's fixed point theorem. The proof  is rather standard but, for the reader's convenience  we give a sketch.  
\begin{Lemme}[Picard's fixed point]\label{Lem-Picard} Let $(E, \Vert \cdot \Vert_E)$ be a Banach space and let $u_0 \in E$ be an initial datum such that $\Vert u_0 \Vert_{E} \leq \delta$. Assume that $B: E \times E \to E$ is a bilinear application and that $L: E \to E$ is a linear application. Assume moreover the following controls:
	\begin{equation}\label{Controls_Picard}
	\Vert B(u,v) \Vert_{E} \leq C_B \Vert u \Vert_{E}\Vert u \Vert_{E}, \qquad \Vert L(u) \Vert_{E} \leq C_L \Vert u \Vert_{E},
	\end{equation}
	for all $u,v \in E$, where the continuity constants of these applications verify:
	\begin{equation}\label{Conditions_Picard}
	0<3  C_L <1, \quad 0< 9 \delta C_B <1 \quad \text{and}\quad  C_L + 6 \delta C_B <1.
	\end{equation}
	Then the equation $u=u_0 + B(u,u)+L(u)$ admits a unique solution $u \in E$ such that $\Vert u \Vert_{E} \leq 3\delta$. 
\end{Lemme}
\pv From the initial datum $u_0$  we consider the iterative equation: $u_{n+1}= u_0 + B(u_n, u_n) + L(u_n)$ for $n \in \mathbb{N}$. By the boundness of the operators $B(\cdot , \cdot)$ and $L(\cdot)$ given in  (\ref{Controls_Picard}), and moreover, by the first and the second  condition in (\ref{Conditions_Picard}), we can prove that $\Vert u_{n+1} \Vert_{E} \leq 3 \delta$ for all $n\in \mathbb{N}$. On the other hand, always by  (\ref{Controls_Picard}) we have the  estimate  $\Vert u_{n+1} - u_{n} \Vert_{E} \leq (C_L + 6 \delta C_B )^n \Vert u_{1} - u_{0} \Vert_{E}.$ Finally, since  $0< C_L + 6 \delta C_B  < 1$  the rest of the proof follows from  well-known arguments.  \finpv 

\medskip

Withing the framework of this lemma, we set the Banach  space
\begin{equation*}
E=\big\{ (\vu, \theta) \in L^{3,\infty}(\Rt): \ \ \text{div}(\vu)=0 \big\},
\end{equation*}
with its natural norm $\|  (\vu, \theta)\|_{L^{3,\infty}}= \| \vu \|_{L^{3,\infty}}+\|\theta\|_{L^{3,\infty}} $. Moreover, from the equations (\ref{Boussineq-Int-1}) and (\ref{Boussineq-Int-2}) we set 
\begin{equation}\label{u-u0}
u= \left( \begin{array}{c}
\vu \\ 
\theta
\end{array}  \right), \quad u_0= \left( \begin{array}{c}\vspace{2mm}
\frac{1}{-\Delta}\Big( \P (\vf) \ \\ 
\frac{1}{-\Delta} (g)
\end{array}  \right),
\end{equation}
and 
\begin{equation}\label{B-L}
B(u,u)= \left( \begin{array}{c}\vspace{2mm}
 - \frac{1}{-\Delta} \left( \P ( \text{div} (\vu \otimes \vu) ) \right)\\ 
- \frac{1}{-\Delta} \left( \text{div}(\theta \vu) \right) 
\end{array}  \right), \quad L(u)= \left( \begin{array}{c}
\frac{1}{-\Delta}\Big( \P \left(  \theta \vec{\bf g} \right)\Big)  \\ 
0
\end{array}  \right).
\end{equation}
We shall construct  a solution to the fixed point problem:
\begin{equation}\label{Fixed-point}
u=u_0 + B(u,u)+ L(u).
\end{equation}

By the hypothesis of Proposition \ref{Prop-Existence-L3}  we directly have $u_0 \in E$, so  we must verify the controls given in (\ref{Controls_Picard}) with constants $C_B, C_L>0$ satisfying (\ref{Conditions_Picard});  and where the quantity $\delta>0$ is defined in (\ref{Control-Forces-Grav}). 

\medskip 

Let us start by studying each component of the bilinear term $B(u,u)$. We remark that in the Fourier level  the operator $\ds{ \frac{\text{div}}{-\Delta}\Big( \P \left(  \cdot \right)\Big)}$ has a symbol $\widehat{m}_1(\xi)=\big(\widehat{m}_{1,i,j,k}(\xi)\big)_{1\leq i,j,k\leq 3}$,  where each $\widehat{m}_{1,i,j,k}$ is a $\mathcal{C}^{\infty}-$ function outside the origin, and moreover, it is a homogeneous function of degree $-1$. Consequently, in the spatial variable, the action of this operator can be seen as a product of convolution  with the tensor $m_1=\big(m_{1,i,j,k}\big)_{1\leq i,j,k\leq 3}$, where  $m_{1,i,j,k}$ is a $\mathcal{C}^{\infty}-$ function outside the origin and homogeneous of degree $-2$.  In particular we have $m_1 \in L^{\frac{3}{2}, \infty}(\Rt)$. 

\medskip

Then, by the first point of Lemma \ref{Lem-Young} (with $p=3$, $p_1=p_2=3/2$ and $q=q_1=q_2=+\infty$) and by the second point of Lemma \ref{Lem-Holder}  (with $p=3/2$, $p_1=p_2=1/3$ and $q=q_1=q_2=+\infty$) we write
\begin{equation*}
\begin{split}
\left\| \frac{1}{-\Delta} \left( \P ( \text{div} (\vu \otimes \vu) ) \right) \right\|_{L^{3,\infty}}  = & \,  \| m_1 \ast (\vu \otimes \vu) \|_{L^{3,\infty}}  \leq \,  C_1 \| m_1 \|_{L^{3/2,\infty}} \| \vu \otimes \vu \|_{L^{3/2,\infty}} \\
\leq &\, C_1 \| m_1 \|_{L^{3/2,\infty}}  \, \| \vu \|^{2}_{L^{3,\infty}}, 
\end{split} \qquad 
\end{equation*}
where we set the constant  $C_{B,1}=C_1 \| m_1 \|_{L^{3/2,\infty}}$. 

\medskip 

The other term $\ds{- \frac{1}{-\Delta} \left( \text{div}(\theta \vu) \right)}$ follows the same estimates: we denote by $m_2 \in L^{3/2,\infty}(\Rt)$ the convolution tensor of the operator $\ds{\frac{1}{-\Delta} \left( \text{div}(\cdot) \right) }$ and we write
\begin{equation*}
\left\|  - \frac{1}{-\Delta} \left( \text{div}(\theta \vu) \right)  \right\|_{L^{3,\infty}} \leq C_1 \| m_2 \|_{L^{3/2,\infty}}\, \| \theta \|_{L^{3,\infty}}\, \| \vu \|_{L^{3,\infty}}, \quad C_{B,2}= C_1 \| m_2 \|_{L^{3/2,\infty}}.
\end{equation*}

With these estimates at hand, we set the constant $C_B= \max(C_{B,1}, C_{B,2})$ and we get 
\begin{equation}\label{Continuity-bilinear}
\| B(u,u) \|_{L^{3,\infty}} \leq C_B  \, \| u \|^{2}_{L^{3,\infty}}. 
\end{equation}

We study now the linear term $L(u)$. In this case we remark that the operator $\ds{\frac{1}{-\Delta}\left( \P (\cdot)\right)}$ has a symbol $\widehat{m}_3(\xi)$ which is a tensor containing   $C^{\infty}$- functions outside the origin and homogeneous of degree $-2$. Then, in the spatial variable we have $\ds{\frac{1}{-\Delta}\left( \P (\theta \vec{{\bf g}})\right)= m_3 \ast  \theta \vec{{\bf g}}}$,  where the tensor $m_3$ is smooth outside the origin and homogeneous of degree $-1$. Consequently, it  verifies $m_3 \in L^{3,\infty}(\Rt)$. 

\medskip 

By the second point of Lemma \ref{Lem-Young} and by the first point of Lemma \ref{Lem-Holder} (with $p_1=3/2$, $p_2=3$, $q_1=1$, $q_2=+\infty$) we write 
\begin{equation*}
\begin{split}
\left\| \frac{1}{-\Delta}\left( \P (\theta \vec{{\bf g}})\right) \right\|_{L^{3,\infty}}= &\, \| m_3 \ast  \theta \vec{{\bf g}} \|_{L^{3,\infty}} \leq C_2 \| m_3 \|_{L^{3,\infty}} \| \theta \vec{{\bf g}} \|_{L^1} \\
\leq &\, C_2  \| m_3 \|_{L^{3,\infty}} \| \vec{{\bf g}} \|_{L^{3/2,1}}\, \| \theta \|_{L^{3,\infty}}.
\end{split}
\end{equation*}
We set the constant $C_L=  C_2 \| m_3 \|_{L^{3,\infty}} \| \vec{{\bf g}} \|_{L^{3/2,1}}$ and we have  the estimate 
\begin{equation}\label{Continuity-Linear}
\| L(u) \|_{L^{3,\infty}} \leq C_L \| u \|_{L^{3,\infty}}. 
\end{equation}

Once we have the constants $C_B$ and $C_L$ given above, we must verify the constraints given in (\ref{Conditions_Picard}). To verify the first constraint we set $\ds{\| \vec{{\bf g}}  \|_{L^{3/2,1}}< \frac{1}{3 C_2 \| m_3 \|_{L^{3,\infty}}}}$ to get $C_L < \frac{1}{3}$. With this inequality the second  and the third constraint in  (\ref{Conditions_Picard})  are  satisfied as long as  $\ds{\delta < \frac{1}{9 C_B}}$. We thus set 
\begin{equation}\label{epsilon1}
\varepsilon_1 < \min\left(  \frac{1}{3 C_2 \| m_3 \|_{L^{3,\infty}}}, \frac{1}{9 C_B}\right),
\end{equation}
 and  Proposition \ref{Prop-Existence-L3} follows from the assumption (\ref{Control-Forces-Grav}) and Lemma \ref{Lem-Picard}. \finpv

\medskip

With this first result, now we are able to prove that the solution  $(\vu,\theta)\in L^{3,\infty}(\Rt)$ also belongs to the space $L^p(\Rt)$ with $3<p<+\infty$, provided that the external forces verify suitable supplementary hypotheses.  In the  case of the Lorentz space $L^{p,q}(\Rt)$ (with $1\leq q \leq +\infty$) the proof follows the same lines, so it is enough to   focus on the Lebesgue spaces. 

\begin{Proposition}\label{Prop-Existence-Lp}  Let $(\vu, \theta)  \in L^{3,\infty}(\Rt)$ be  the solution to  the equations (\ref{Boussineq-Int-1})-(\ref{Boussineq-Int-2}) constructed in  Proposition \ref{Prop-Existence-L3}. Let $\delta>0$ be the quantity given in (\ref{Control-Forces-Grav}) and $\varepsilon_1>0$ defined in (\ref{epsilon1}). 
	
\medskip 	
	
Assume in addition that
\begin{equation}\label{Forces-Lp}
\begin{cases}\vspace{2mm}
\frac{1}{-\Delta}\Big( \P(\vf) \Big) \in L^{p} \cap L^4 (\Rt), \ \  \frac{1}{-\Delta}(g) \in L^{p} \cap L^4 (\Rt), \ \ 3<p<4, \\
\frac{1}{-\Delta}\Big( \P(\vf) \Big) \in L^{p}  (\Rt), \ \  \frac{1}{-\Delta}(g) \in L^{p} (\Rt), \ \ 4 \leq p \leq +\infty, 
\end{cases}
 \ \ \ \mbox{and} \ \ \vec{{\bf g}} \in L^{3,1}(\Rt). 
\end{equation}
 There exists a universal quantity $0<\varepsilon<\varepsilon_1$ (defined in (\ref{epsilon})), such that if $\delta <\varepsilon$ then we have  $(\vu, \theta)  \in L^p(\Rt)$  with $3<p\leq +\infty$. 
\end{Proposition}
\pv  Recall that the solution $(\vu, \theta)  \in L^{3,\infty}(\Rt)$ to the equations (\ref{Boussineq-Int-1})-(\ref{Boussineq-Int-2}) is obtained as the limit (in the strong topology of the space $L^{3,\infty}(\Rt)$) of the sequence $(u_n)_{n\in \mathbb{N}}$, which  defined by  the iterative expression: 
\begin{equation}\label{Sequence}
u_{n+1}= u_0+B(u_n, u_n)+L(u_n), \qquad \mbox{for} \ \ n \geq 0,
\end{equation}
 and where  $u_0$, $B(u_n, u_n)$ and $L(u_n)$ are given in the expressions (\ref{u-u0}) and (\ref{B-L}) respectively.  We shall use this sequence to prove that $(\vu, \theta) \in L^p(\Rt)$.  We thus write
 \begin{equation}\label{Sequence-Lp}
 \| u_{n+1} \|_{L^p} \leq \| u_0 \|_{L^p} + \| B(u_n u_n) \|_{L^p} + \| L(u_n) \|_{L^p},
 \end{equation}
 where, by the assumption (\ref{Forces-Lp}) we have $\| u_0 \|_{L^p}  < +\infty $. On the other hand, in order to estimate the bilinear term we recall that its first component writes down as $m_1 \ast (\vu \otimes \vu)$ with $m_1 \in L^{3/2,\infty}(\Rt)$. Then, by the first point of Lemma \ref{Lem-Young} (with $p=q$, $p_1=3/2$, $p_2=3p/(3+p)$, $q_1=+\infty$ and $q_2=p_2$)  we get
 \begin{equation}\label{C1}
 \| m_1 \ast (\vu \otimes \vu)\|_{L^p} \leq C_1(p) \,  \| m_1 \|_{L^{3/2,\infty}} \| \vu \otimes \vu \|_{L^{\frac{3p}{3+p}}}, \quad \mbox{with} \ \  0<C_1(p)=C\frac{3p^2}{2p-3}.
 \end{equation}
 \begin{Remarque}\label{Rmk-C1} Note that the constant $C_1(p)$ given above blows-up at $p=\frac{3}{2}$ and $p=+\infty$. This fact yields the first  constraint $\frac{3}{2}<p<+\infty$, which is the same in \cite{Bjorland}.
 \end{Remarque}	
Then, from the estimate (\ref{C1}) and by second point of Lemma \ref{Lem-Holder} (remark that $(3+p)/3p=1/3+1/p$) we write
\begin{equation*}
\| m_1 \ast (\vu \otimes \vu) \|_{L^p} \leq C_1(p) \| m_1 \|_{L^{3/2,\infty}}\,  \| \vu \|_{L^p}\, \| \vu \|_{L^{3,\infty}}. 
\end{equation*}
Moreover, recall that the sequence $(u_n)_{n \in \mathbb{N}}$ defined in (\ref{Sequence}) verifies $\| u_n \|_{L^{3,\infty}} \leq 3\delta$ and  we thus obtain 
\begin{equation*}
\| m_1\ast (\vu \otimes \vu) \|_{L^p} \leq 3\delta C_1(p)\, \| m_1 \|_{L^{3/2,\infty}}\, \| \vu \|_{L^p}.  
\end{equation*} 

The second component of the bilinear form $B(u_n, u_n)$ essentially follows the same arguments and we have
\begin{equation}\label{Estim-B2}
\left\|  \frac{1}{-\Delta} \text{div} (\theta \vec{\bf g}) \right\|_{L^{p}}=  \| m_2\ast (\theta  \vu) \|_{L^p} \leq 3\delta C_1(p)\, \| m_2 \|_{L^{3/2,\infty}}\, \| \theta \|_{L^p}.   
\end{equation} 

With these estimates at our disposal, for the sake of simplicity we introduce the constant  $C_m=\max(\| m_1 \|_{L^{3/2,\infty}},  \| m_2 \|_{L^{3/2,\infty}})$
 and we obtain the estimate 
\begin{equation}\label{Bilinear-Lp}
\| B(u_n,u_n) \|_{L^p} \leq 3\delta C_m  C_1(p)  \| u_n \|_{L^p}. 
\end{equation}

We study now the linear term $L(u_n)$. Let $1<p_1<3/2$ defined by the relationship $1/p_1= 2/3+1/p$. 
\begin{Remarque}\label{Rmk-2} Observe that the relationship $1/p_1= 2/3+1/p$ yields the new  constraint $3<p<+\infty$.
\end{Remarque}	
By the first point of Lemma \ref{Lem-Young} (since  $1/p_1= 2/3+1/p$ we have $1+1/p=1/3+1/p_1$, hence we set  $p_2=3$, $q_2=+\infty$ and  $q=p$, $q_1=p_1$)  and by using again the second point of  Lemma \ref{Lem-Holder}  (recall that $1/p_1=2/3+1/p$) we write
\begin{equation*}
\begin{split}
\left\| \frac{1}{-\Delta} \left( \P(\theta \vec{{\bf g}}) \right) \right\|_{L^p}=&\, \| m_3 \ast (\theta \vu) \|_{L^p} \leq C{'}_{1}(p) \| m_3 \|_{L^{3,\infty}}\,  \| \theta \vec{{\bf g}} \|_{L^{p_1}} \leq C^{'}_{1}(p) \| m_3 \|_{L^{3,\infty}}\, \| \vec{{\bf g}} \|_{L^{3/2}}\, \| \theta \|_{L^p} \\
\leq &\,  C^{'}_{1}(p) \| m_3 \|_{L^{3,\infty}}\, \| \vec{{\bf g}} \|_{L^{3/2,1}}\, \| \theta \|_{L^p}.
\end{split}
\end{equation*}
\begin{Remarque}\label{Rmk3}
We have $C^{'}_{1}(p)= C\frac{3p^2}{p-3}$,  which is well-defined as long as  $3<p<+\infty$.
\end{Remarque}
For the sake of simplicity, we introduce the constant $C^{'}_{m}=\| m_3 \|_{L^{3,\infty}}$ and we thus  get
\begin{equation}\label{Linear-Lp}
\| L(u_n) \|_{L^p} \leq   \left( C^{'}_{m}C^{'}_{1}(p)   \| \vec{{\bf g}} \|_{L^{3/2,1}}\right) \, \| u_n \|_{L^p}. 
\end{equation}

Once we have the estimates (\ref{Bilinear-Lp}) and (\ref{Linear-Lp}), we get back to the inequality (\ref{Sequence-Lp}) to write
\begin{equation}\label{Sequence-Lp-2}
\| u_{n+1}\|_{L^p} \leq \|u_0 \|_{L^p}+ \left(  3\delta C_m  C_1(p)  + C^{'}_{m}C^{'}_{1}(p) \, \| \vec{{\bf g}} \|_{L^{3/2,1}}\right) \| u_n \|_{L^p}, \quad \mbox{for} \ \ n \geq 0 \ \ \mbox{and} \ \ 3<p<+\infty. 
\end{equation}

At this point, we need to find  additional constraints on the quantities  $\delta$  and $\| \vec{{\bf g}} \|_{L^{3/2,1}}$  to obtain 
\begin{equation}\label{Inequality}
\left(  3\delta C_m  C_1(p)  + C^{'}_{m}C^{'}_{1}(p) \, \| \vec{{\bf g}} \|_{L^{3/2,1}}\right) < \frac{1}{2}.
\end{equation}
With this inequality we will able to prove that the sequence $(u_n)_{n\in \mathbb{N}}$ is uniformly bounded in the space $L^p(\Rt)$ (see the estimate (\ref{Uniform-control-un}) below) and consequently its limit verifies $u\in L^p(\Rt)$.  However, there is a problem to overcome: the quantities $C_1(p)$ and  $C^{'}_{1}(p)$ (defined in the expression (\ref{C1}) and the Remark \ref{Rmk3} respectively) are not bounded in the whole interval $3<p<+\infty$.  So, we shall use the following strategy: first we will prove the wished inequality (\ref{Inequality}) in the interval $4\leq p \leq 7$ to obtain $u \in L^p(\Rt)$ for these values of $p$. Thereafter, we shall use this information to get that $u \in L^p(\Rt)$ in the intervals $3<p<4$, $7<p<+\infty$ and the value $p=+\infty$. 

\begin{itemize}
	\item {\bf The case $4\leq p \leq 7$}.   We get back to the expression of the quantity $C_1(p)$ given in (\ref{C1}) and we define  $\ds{0<M=\max_{ 4 \leq  p\in 7} C_1(p)<+\infty}$. Then, we set the new constraint on the parameter $\delta$ (which already verifies (\ref{Control-Forces-Grav})) as follows:
	\begin{equation}\label{Condition-delta-2}
	3  \delta C_m M  < \frac{1}{4}.  
	\end{equation}
	Similarly, we  get back to the expression of the quantity $C^{'}_{1}(p)$ given in Remark \ref{Rmk3} and we define  $\ds{0<M'=\max_{4\leq  p \leq 7} C^{'}_1(p)<+\infty}$. We thus  set the new constraint on the quantity $\| \vec{{\bf g}} \|_{L^{3/2,1}}$ (which also verifies (\ref{Control-Forces-Grav})) as the next one:
	\begin{equation}\label{Condition-g-2}
	C^{'}_{m} M' \| \vec{{\bf g}} \|_{L^{3/2,1}} < \frac{1}{4}.
	\end{equation}
By these new constraints,  for all $4 \leq p \leq 7$ we have the inequality (\ref{Inequality});   and we get back to the inequality (\ref{Sequence-Lp-2})  to write 
	\begin{equation*}
	\| u_{n+1} \|_{L^p} \leq  \| u_0 \|_{L^p} + \frac{1}{2} \| u_{n} \|_{L^p},  \ \ \mbox{for all} \ \ n\geq 0.   
	\end{equation*}
	Hence, we iterate to  obtain the uniform control
	\begin{equation}\label{Uniform-control-un}
	\| u_{n+1} \|_{L^p} \leq \left( \sum_{k=0}^{+\infty} \frac{1}{2^k} \right) \, \| u_0 \|_{L^p}, \ \ \mbox{for all} \ \ n\geq 0,  
	\end{equation}
and then we have $ u=(\vu,\theta) \in L^p(\Rt)$ with $4 \leq p \leq 7$. 

\medskip 

 Moreover, we emphasize that this fact holds true as long as the  quantities $\delta$ and $ \| \vec{{\bf g}} \|_{L^{3/2,1}} $ verify the set of constraints  (\ref{Control-Forces-Grav}), (\ref{Condition-delta-2}) and (\ref{Condition-g-2}), which can be jointly written as
	\begin{equation}\label{epsilon} 
	\delta,   \| \vec{{\bf g}} \|_{L^{3/2,1}} <  \min\left( \frac{1}{12 C_m M} ,  \frac{1}{4 C^{'}_{m} M'}, \varepsilon_1 \right)=\varepsilon.     
	\end{equation}
\item {\bf The case $3<p<4$}. By the  assumption (\ref{Forces-Lp}) we have $u_0 \in L^4(\Rt)$ and by the uniform control (\ref{Uniform-control-un}) we have $u\in L^4(\Rt)$. Moreover, recall that by Proposition (\ref{Prop-Existence-L3}) we also have $u \in L^{3,\infty}(\Rt)$ and by the second point of Lemma \ref{Lem-Interpolation} we have $u=(\vu,\theta)\in L^p(\Rt)$ with $3<p<4$.  	
	
\item {\bf The case $ 7 <p < +\infty $}. First recall that by assumption (\ref{Forces-Lp}) we have $u_0 \in L^p(\Rt)$ and since we also $u_0 \in L^{3,\infty}(\Rt)$ by the standard interpolation inequalities  we obtain $u_0 \in L^{4} \cap  L^{7}(\Rt)$. Thus, always by the uniform control (\ref{Uniform-control-un}) we obtain $u \in L^{4} \cap L^{7}(\Rt)$  and then by the second point of Lemma \ref{Lem-Interpolation} we have $u \in L^{6,2}(\Rt)$.    

\medskip

With this information we can prove that $B(u, u) \in L^\infty(\Rt)$. Indeed, this bilinear form essentially writes down as $m\ast (u \otimes u)$  where the tensor $m\in L^{3/2,\infty}(\Rt)$ denotes each tensor $m_i$ (with $i=1,2$) in each component of $B(u,u)$ (see the expression (\ref{B-L})). Then, since $u \in L^{6,2}(\Rt)$ by the second point of Lemma \ref{Lem-Holder} we obtain $u \otimes u \in L^{3,1}(\Rt)$; and since $m \in L^{3/2,\infty}(\Rt)$ by the third point of  Lemma \ref{Lem-Young}  we have 
\begin{equation}\label{Bilinear-Linfty}
B(u, u) \in L^{\infty}(\Rt).
\end{equation}
As we also have $B(u, u) \in L^{3,\infty}(\Rt)$ (see the estimate (\ref{Continuity-bilinear})) and by using again the second point of Lemma \ref{Lem-Interpolation} we obtain $B(u,u) \in L^p(\Rt)$ with $p \in (7,+\infty)$. 

\medskip

Now, we shall verify that $L(u) \in L^{\infty}(\Rt)$. Indeed, by the third point of Lemma \ref{Lem-Young}  we write
\begin{equation*}
\left\| \frac{1}{-\Delta} \left( \P (\theta \vec{{\bf g}})\right) \right\|_{L^\infty}=\left\| m_3\ast  (\theta \vec{{\bf g}})  \right\|_{L^\infty} \leq C_3 \, \| m_3 \|_{L^{3,\infty}}\, \| \theta \vec{{\bf g}} \|_{L^{3/2,1}}.
\end{equation*}
But, since $\theta \in L^{3,\infty}(\Rt)$  and   $\vec{{\bf g}} \in L^{3,1}(\Rt)$  (see the assumption (\ref{Forces-Lp})) by the second point of Lemma \ref{Lem-Holder}  we get $\| \theta \vec{{\bf g}} \|_{L^{3/2,1}} \leq  \| \theta \|_{L^{3,\infty}}\| \vec{{\bf g}} \|_{L^{3,1}}$. We thus obtain
\begin{equation}\label{Linear-Linfty}
L(u) \in L^{\infty}(\Rt).
\end{equation}
As before, since we also have $L(u)\in L^{3,\infty}(\Rt)$ (see the estimate (\ref{Continuity-Linear})) we use  using the  second point of Lemma \ref{Lem-Interpolation} to conclude that  $L(u)\in L^p(\Rt)$ with $p\in (7,+\infty)$. Consequently, we have $u=(\vu,\theta)\in L^p(\Rt)$ with $7<p<+\infty$. 

\medskip

\item {\bf The case $p=+\infty$}. We just remark that from the information given in (\ref{Bilinear-Linfty}) and (\ref{Linear-Linfty}) we have $u=(\vu,\theta) \in L^{\infty}(\Rt)$ as long as $u_0 \in L^\infty(\Rt)$, which holds true under the assumption $\ds{\frac{1}{-\Delta}\Big(\P(\vf)\Big) \in L^{\infty} (\Rt)}$ and    $\ds{\frac{1}{-\Delta}(g) \in L^{\infty} (\Rt)}$.  Proposition \ref{Prop-Existence-Lp} is proven.  \finpv 
\end{itemize}	 
Once we have  Propositions \ref{Prop-Existence-L3} and \ref{Prop-Existence-Lp},  Theorem \ref{Th1} is now proven. \finpv 

\section{Proof of Proposition \ref{Prop-Nonexistence}}\label{Sec-Non-existence} 
The proof is divided into three main steps,  which we shall  study separately  in the next (technical) propositions. In addition, we recall that in all this section the external forces $\vf, g \in \mathcal{S}(\Rt)$ and the gravitational acceleration $\vec{\bf g} \in \mathcal{S}(\Rt)$ verify (\ref{Force-Supp-Fourier}) and we have $\ds{(-\Delta)^{-1} \vf\in \mathcal{S}(\Rt)},   \ds{(-\Delta)^{-1} g \in \mathcal{S}(\Rt)}$ and $ (-\Delta)^{-1} \vec{{\bf g}} \in \mathcal{S}(\Rt)$. 

\medskip

{\bf Step 1. Pointwise decaying properties}. We start by proving the following pointwise decaying properties of solutions to the system (\ref{Boussinesq-toy-model}). 
\begin{Proposition}\label{Prop-Pointwise}  Let $\vf, g, \vec{\bf g} \in \mathcal{S}(\Rt)$  verifying  (\ref{Force-Supp-Fourier}).   Moreover, let  $\varepsilon>0$ be the quantity defined in (\ref{epsilon}). 

\medskip

There exists a universal quantity $0<\varepsilon_2  < \varepsilon$, such that if 
\begin{equation}\label{delta2}
\delta_2= \| (-\Delta)^{-1} \vf \|_{L^{3,\infty}}+ \| \, |x|  (-\Delta)^{-1}  \vf \|_{L^{\infty}}+  \| (-\Delta)^{-1} g \|_{L^{3,\infty}}+ \| \, |x|  (-\Delta)^{-1} g \|_{L^{\infty}} < \varepsilon_{2}, 
\end{equation}	
and 
\begin{equation}\label{Cond-g-2}
\| \vec{{\bf g}} \|_{L^{3/2,1}}+\| (-\Delta)^{-1} | \vec{{\bf g}}| \|_{L^{\infty}} < \varepsilon_{2},
\end{equation}	
then the solution  $(\vu, \theta)\in L^{3,\infty}(\Rt)\cap L^\infty(\Rt)$  to the system (\ref{Boussinesq-toy-model}), given by Theorem \ref{Th1},   verifies the  following pointwise estimates:
\begin{equation}\label{Pointwise}
|\vu(x)| \leq \frac{C_1}{1+|x|}, \quad  |\theta(x)| \leq \frac{C_2}{1+|x|}, \qquad x \in \Rt,
\end{equation}	
with constants $C_1, C_2>0$ depending on $\vf, g, \vec{{\bf g}}, \vu$ and $\theta$, but independent on $x$. 
\end{Proposition}
\begin{Remarque}  It is worth emphasizing this result also holds true for the (original) system (\ref{Boussinesq}), since here we only need the information $\varphi \in L^\infty(\Rt)$; and we thus can set $\varphi \equiv 1$. 
\end{Remarque}	
\pv First,  remark  that the solution $(\vu, \theta) \in L^{3,\infty}(\Rt)$ is   constructed as in the proof of Proposition \ref{Prop-Existence-L3}:  we assume the stronger inequality (\ref{delta2}) (remark that $\delta < \delta_2$, where $\delta$ is defined in (\ref{Control-Forces-Grav}))  and we solve the fixed point equation 
\begin{equation}\label{Fixed-point-varphi}
u = u_0 + B_\varphi (u,u)+ L(u)
\end{equation}
in the Banach space  $L^{3,\infty}(\Rt)$. Here,  the terms  $u, u_0$ are given in (\ref{u-u0}) and the term  $L(u)$ is given in  (\ref{B-L}) (without the Leray's projector $\P$). Moreover we have 
\begin{equation}\label{Bvarphi}
B_\varphi (u,u)=\left( \begin{array}{c}\vspace{2mm}
- \frac{1}{-\Delta} \left( \P ( \text{div} (\varphi \vu \otimes \vu) ) \right)\\ 
- \frac{1}{-\Delta} \left( \text{div}(\theta \vu) \right) 
\end{array}  \right). 
\end{equation}
This solution is uniquely defined by the condition $\| u \|_{L^{3,\infty}} \leq 3 \delta_2$. 

\medskip 

On the other hand, we introduce the weighted $L^{\infty}$- space: 
\begin{equation}
L^{\infty}_{1}(\Rt)=\big\{  f : \Rt \to \Rt \ / \ \ f \ \mbox{is measurable and} \ \  |x| f \in L^{\infty}(\Rt) \big\},
\end{equation}
which verifies $L^{\infty}_{1}(\Rt) \subset L^{3,\infty}(\Rt)$. Thus, we shall solve the problem (\ref{Fixed-point-varphi})  in the smaller Banach space $L^{\infty}_{1}(\Rt)$ with the norm 
\begin{equation*}
\|  f \|_{1}= \| \,  |x| f \|_{L^\infty}+\| f \|_{L^{3,\infty}}. 
\end{equation*}

Since $\varphi \in L^{\infty}(\Rt)$   the quantity $\| B_\varphi (u,u)\|_{L^{3,\infty}}$ is similarly estimated as in (\ref{Continuity-bilinear}) and we have
\begin{equation}\label{Continuity-Bilinear-varphi1}
\| B_\varphi (u,u) \|_{L^{3,\infty}} \leq \| \varphi \|_{L^\infty} C_B \| u \|^{2}_{L^{3,\infty}}. 
\end{equation}
Moreover, the quantity $\| L(u) \|_{L^{3,\infty}}$ was already estimated in (\ref{Continuity-Linear}). So, it remains to estimate the quantities $\| \, |x| B_\varphi(u,u) \|_{L^{\infty}}$ and $\|\, |x| L(u)\|_{L^\infty}$. The first quantity follows the same estimates in \cite[Lemma $3.6$]{Bjorland} and we can write 
\begin{equation}\label{Continuity-Bilinear-varphi2}
\| \, |x| B_\varphi (u,u) \|_{L^\infty} \leq C\| \varphi \|_{L^\infty} \| \, |x| u \|^{2}_{L^\infty}. 
\end{equation}
To study the second quantity, for $x\neq 0$ fixed  we write 
\begin{equation*}
\begin{split}
\left| \frac{1}{-\Delta} (\theta \vec{\bf {g}}) (x) \right| \leq &\, \int_{\Rt} \frac{1}{|x-y |} | \theta (y) |\, |\vec{{\bf g}}(y) | dy \leq  \left(\int_{|y|\leq |x|/2} + \int_{|y| >  |x|/2} \right)  \frac{1}{|x-y|} | \theta (y) |\, |\vec{{\bf g}}(y) | dy \\
=&\, I_1 + I_2.
\end{split}
\end{equation*}
To estimate the term $I_1$, we remark that since $| y | \leq |x|/2$ we have $|x-y|\geq |x|/2$ and by the first point of Lemma \ref{Lem-Holder} we can write 
\begin{equation*}
I_1 \leq \frac{C}{|x|} \int_{|y|\leq |x|/2} | \theta (y) |\, |\vec{{\bf g}}(y) | dy  \leq \|  \theta \vec{{\bf g}} \|_{L^1} \leq \frac{C}{|x|}\|\theta \|_{L^{3,\infty}}\,\| \vec{{\bf g}} \|_{L^{3/2,1}}. 
\end{equation*}
To estimate the term $I_2$, since $| y | > |x|/2$ we can use  expression $\| \theta \|_{L^{\infty}_{1}}$ to get
\begin{equation*}
\begin{split}
I_2 \leq &\,  C \| \theta \|_{L^{\infty}_{1}} \int_{| y | > |x|/2} \frac{1}{|x-y|\,|y|} |\vec{{\bf g}}(y) | dy  \leq  \frac{C}{|x|} \| \theta \|_{L^{\infty}_{1}} \int_{| y | > |x|/2} \frac{1}{|x-y|} |\vec{{\bf g}}(y) | dy \\
\leq & \, \frac{C}{|x|} \| \theta \|_{L^{\infty}_{1}} \int_{\Rt} \frac{1}{|x-y|} |\vec{{\bf g}}(y) | dy \leq \frac{C}{|x|} \| \theta \|_{L^{\infty}_{1}} (-\Delta)^{-1}| \vec{{\bf g}}(x) | \leq \frac{C}{|x|} \| \theta \|_{L^{\infty}_{1}} \| (-\Delta)^{-1}| \vec{{\bf g}}| \|_{L^{\infty}}. 
\end{split}
\end{equation*}
From these estimate we obtain
\begin{equation}\label{Continuity-Linear-2}
\|\, |x| L(u) \|_{L^{\infty}} \leq C\, (\| \vec{{\bf g}} \|_{L^{3/2,1}}+ \| (-\Delta)^{-1}| \vec{{\bf g}}| \|_{L^{\infty}}) \| \theta \|_{1}.  
\end{equation}

Once we have the estimates (\ref{Continuity-Bilinear-varphi1}), (\ref{Continuity-Bilinear-varphi1})  and (\ref{Continuity-Linear-2}), and moreover, by the estimate (\ref{Continuity-Linear}) we are able to write
\begin{equation*}
\| B_\varphi (u,u) \|_1 \leq \max(C_B, C)\| \varphi \|_{L^\infty} \| u \|^{2}_{1}, \quad \| L(u) \| \leq \max(C,C_L) (\| \vec{{\bf g}} \|_{L^{3/2,1}}+ \| (-\Delta)^{-1}| \vec{{\bf g}}| \|_{L^{\infty}}) \| u \|_{1}.
\end{equation*}

We  set a quantity $0<\varepsilon_{2}<\varepsilon$ small enough such that the constraints given in (\ref{Conditions_Picard}) hold. We also assume (\ref{delta2} ) and (\ref{Cond-g-2}).  Consequently, by Lemma  \ref{Lem-Picard} there exists a solution $u=(\vu, \theta) \in L^{\infty}_{1}$ to the equation (\ref{Fixed-point-varphi})  which is  uniquely determined by the condition $\| u \|_{1}\leq 3 \delta_2$. But, since $\| u \|_{L^{3,\infty}} \leq \| u \|_1$ this solution is the same to the one given by Theorem \ref{Th1}.

\medskip

Finally, recall that this solution also verifies $u \in L^{\infty}(\Rt)$, hence we obtain the wished estimates  (\ref{Pointwise}). Proposition \ref{Prop-Pointwise} is proven.  \finpv  

\medskip

{\bf Step 2. Asymptotic profile}.  Once we have the estimates (\ref{Pointwise}), in the next proposition we develop an asymptotic profile for the velocity $\vu$. 

\begin{Proposition}\label{Prop-Profile} Within the framework of Proposition \ref{Prop-Pointwise}, the velocity $\vu$ has the asymptotic behavior
\begin{equation}\label{Asymptotic-Profile}
\vu(x)= \frac{1}{|x|} \left( \int_{\Rt} \theta(y) \vec{{\bf g}}(y) dy \right)+\textit{o}\left(\frac{1}{|x|}\right), \quad |x|\to +\infty. 
\end{equation}	
\end{Proposition}
\begin{Remarque} To the best of our knowledge, the identity (\ref{Asymptotic-Profile}) strongly depends on the well-prepared function $\varphi$. In particular, to justify our estimates, here  we  need  the pointwise inequality $| \varphi (x) | \lesssim \frac{1}{1+|x|^2}$.  
\end{Remarque}	
\pv Recall that $\vu$ verifies the fixed point equation
\begin{equation}
\vu(x)=\frac{1}{-\Delta}(\theta \vec{{\bf g}})(x)- \frac{1}{-\Delta}(\text{div}(\varphi \vu \otimes \vu))(x)+ \frac{1}{-\Delta}(\vf)(x). 
\end{equation}	 
Our starting point is to prove that the first term on the right-hand side splits as follows:
\begin{equation}\label{Iden-o-1}
\frac{1}{-\Delta}(\theta \vec{{\bf g}})(x) = \frac{1}{|x|} \left( \int_{\Rt} \theta(y) \vec{{\bf g}}(y) dy \right)+ \textit{o}\left(\frac{1}{|x|}\right), \quad |x|\to +\infty.
\end{equation}
Indeed, for $|x|>0$ fixed large enough we write 
\begin{equation*}
\begin{split}
\frac{1}{-\Delta}(\theta \vec{{\bf g}})(x) =  &\, \int_{\Rt} \frac{1}{|x-y|} \theta(y) \vec{{\bf g}}(y) dy\\
=&\, \frac{1}{|x|} \left( \int_{\Rt} \theta(y) \vec{{\bf g}}(y) dy \right) + \int_{\Rt} \left(  \frac{1}{|x-y|}  -  \frac{1}{|x|} \right) \theta(y) \vec{{\bf g}}(y) dy. 
\end{split}
\end{equation*} 
The, for $i=1,2,3$ we shall prove that $\ds{\int_{\Rt} \left(  \frac{1}{|x-y|}  -  \frac{1}{|x|} \right) \theta(y) {\bf g}_i(y) dy = \textit{o}\left( \frac{1}{|x|}\right)}$.  We thus write
\begin{equation*}
\int_{\Rt} \left(  \frac{1}{|x-y|}  -  \frac{1}{|x|} \right) \theta(y) {\bf g}_i(y) dy = \left( \int_{|y| \leq |x|/2 } + \int_{|y|>|x|/2} \right) \left(  \frac{1}{|x-y|}  -  \frac{1}{|x|} \right) \theta(y) {\bf g}_i(y) dy = I_1+ I_2.
\end{equation*}
To estimate the term $I_1$, recall that the function $\frac{1}{|x|}$ is smooth outside the origin and by the mean value theorem, for $w=\alpha(x-y)+(1-\alpha)x $ with $0<\alpha<1$,  we write  $\left|  \frac{1}{|x-y|}  -  \frac{1}{|x|}   \right|  \leq \frac{1}{| w |^2} |y|$. Moreover, remark that   $w=x-\alpha y$  and as we have  $| y |\leq |x|/2$  then we  get $|w| \geq |x | - \alpha |y| \geq |x|-|y| \geq |x|/2$. We thus obtain   $\left|  \frac{1}{|x-y|}  -  \frac{1}{|x|}   \right|  \leq \frac{C}{|x|^2} |y|$. With this inequality at hand, and  moreover, by recalling that ${\bf g}_i \in \mathcal{S}(\Rt)$  and $\theta \in L^{\infty}(\Rt)$,  the term $I_1$ verifies
\begin{equation}\label{Estim-o-1}
I_1 \leq \frac{C}{|x|^2} \| \theta \|_{L^\infty}\, \int_{|y|\leq |x|/2} | y | | {\bf g}_i (y) | dy \leq   \frac{C}{|x|^2} \| \theta \|_{L^\infty}\, \| \, | \cdot | {\bf g}_i \|_{L^1}. 
\end{equation}
To estimate the term $I_2$ we write 
\begin{equation*}
I_2 \leq \int_{|y|>|x|/2} \frac{1}{|x-y|} | \theta(y)||{\bf g}_i(y)| dy + \frac{1}{|x|} \int_{|y|>|x|/2} |\theta(y)||{\bf g}_i(y)| dy=I_{2,1}+I_{2,2}.
\end{equation*}
Then, in order to study the term $I_{2,1}$, we recall that $\theta$ verifies the pointwise estimate (\ref{Pointwise}), and moreover, since ${\bf g}_i \in \mathcal{S}(\Rt)$ we have $| {\bf g}_i (y)| \leq \frac{C}{|y|^\beta}$, with $2<\beta <3$. We thus get
\begin{equation}\label{Estim-o-2}
I_{2,1} \leq C \int_{|y|>|x|/2} \frac{dy}{|x-y| |y| |y|^\beta} \leq \frac{C}{|x|} \int_{\Rt} \frac{dy}{|x-y| |y|^\beta} \leq \frac{C}{|x|^{\beta-1}}.  
\end{equation}
Thereafter, in order to estimate the term $I_{2,2}$, we use again the previous arguments to write 
\begin{equation}\label{Estim-o-3}
I_{2,2} \leq \frac{C}{|x|^2} \int_{|y|>|x|/2} |{\bf g}_i (y) | dy \leq \frac{C}{|x|^2} \| {\bf g}_i \|_{L^1}. 
\end{equation}
With the estimates (\ref{Estim-o-1}), (\ref{Estim-o-2}) and (\ref{Estim-o-3}) at hand, we finally obtain the wished identity (\ref{Iden-o-1}). 

\medskip

Now, we need to verify that 
\begin{equation}\label{Iden-o-2}
\frac{1}{-\Delta}(\text{div}(\varphi \vu \otimes \vu))(x)=\textit{o}\left( \frac{1}{|x|}\right), \qquad |x|\to +\infty.
\end{equation}
Since $| \varphi(x) | \leq \frac{C}{1+|x|^2}$ one can follow the same estimates in the proof of \cite[Theorem $3.1$, equation $(2.3)$]{Bjorland}  to obtain the asymptotic profile
\begin{equation*}
\frac{1}{-\Delta}(\text{div}(\varphi \vu \otimes \vu))(x)=m_2(x) \left( \int_{\Rt} \varphi (y) \vu \otimes \vu (y)dy\right)+ \textit{O}\left( \frac{\log(|x|)}{|x|^3} \right), \quad |x|\to +\infty,
\end{equation*}
where $m_2$ denotes the tensor of the operator $\frac{\text{div}(\cdot)}{-\Delta}$; and it verifies the estimate  $| m_2 (x) | \leq \frac{C}{|x|^2}$ outside the origin. Hence we have the identity (\ref{Iden-o-2}).

\medskip

Finally, recall that the condition (\ref{Force-Supp-Fourier}) yields  $\ds{ \frac{1}{-\Delta}(\vf) \in \mathcal{S}(\Rt)}$. With this information and by the identities (\ref{Iden-o-1}) and (\ref{Iden-o-2}) we obtain  the wished asymptotic profile  (\ref{Asymptotic-Profile}). Proposition \ref{Prop-Profile} is proven. \finpv 

\medskip

{\bf Step 3. Estimate from below}. In out last step, the asymptotic profile (\ref{Asymptotic-Profile}) will allow us to obtain the next estimate from below on the velocity $\vu$, which yields  $\vu  \notin L^p(\Rt)$ for $1\leq p \leq 3$.   
\begin{Proposition}\label{Prop-Estimate-Below} Within the framework of Proposition \ref{Prop-Pointwise}, assume that 
\begin{equation}\label{New-force}
\vf = \kappa \vec{\mathfrak{f}}, \quad g = \kappa \mathfrak{g} \quad \mbox{and} \quad \vec{\bf g}=\kappa \vec{{\mathfrak{ g}}},
\end{equation}	 
with   $\vec{\mathfrak{f}}, \mathfrak{g}, \vec{{\mathfrak{ g}}} \in \mathcal{S}(\Rt)$ verifying (\ref{Force-Supp-Fourier}); and with $\kappa>0$  a suitable quantity  defined in  the expression (\ref{kappa}) below. Moreover, for the external force $\mathfrak{g}$ and  gravitational acceleration $\vec{{\mathfrak{ g}}}=(\mathfrak{g}_1, \mathfrak{g}_2, \mathfrak{g}_3)$, assume that 
\begin{equation}\label{Cond-g}
\int_{\Rt}  \frac{1}{-\Delta} (\mathfrak{g})(x) \mathfrak{g}_i(x) dx \neq 0, \qquad i=1,2,3. 
\end{equation}
Then there exists two quantities $M_1=M_1(\theta, \vec{\bf g})>0$ and  $M_2= M(\theta, \vec{\bf g})>0$, which depend on $\theta$ and $g$, such that  the velocity $\vu$ verifies the estimate from below:
\begin{equation}\label{Estim-from-above}
\frac{M_1}{2|x|} \leq | \vu(x)|,  \qquad M_2 < |x|.
\end{equation} 
Then, $\vu \notin L^p(\Rt)$ for $1\leq p \leq 3$. 
\end{Proposition}	
\pv From the asymptotic profile (\ref{Asymptotic-Profile}), for $|x|>0$ large enough  we can write
\begin{equation}\label{Estim-below}
| \vu(x)| = \left| \frac{1}{|x|} \left( \int_{\Rt} \theta(y) \vec{{\bf g}}(y) dy \right)+\textit{o}\left(\frac{1}{|x|}\right)  \right| \geq \frac{1}{|x|} \left| \int_{\Rt} \theta(y) \vec{{\bf g}}(y) dy \right| - \left| \textit{o}\left(\frac{1}{|x|}\right)  \right|. 
\end{equation}
We thus  set the quantity $\ds{M_1=\left| \int_{\Rt} \theta(y) \vec{{\bf g}}(y) dy \right|}$ and, in order to verify that $M_1>0$, we have the following technical lemma:
\begin{Lemme} The assumptions (\ref{New-force}) and  (\ref{Cond-g}) imply that $M_1>0$.
\end{Lemme}
\pv  Recall that the temperature $\theta$ verifies the fixed point equation
\begin{equation*}
\theta = \frac{1}{-\Delta} (g) - \frac{1}{-\Delta}(\text{div}(\theta \vu)),
\end{equation*}
where we  denote $\ds{\frac{1}{-\Delta} (g)=\theta_0}$ and $\ds{\frac{1}{-\Delta}(\text{div}(\theta \vu))=B_2(\theta, \vu)}$.  Moreover, since $g = \kappa \mathfrak{g}$ we shall denote   $\ds{\theta_0 = \kappa \tilde{\theta}_0= \kappa\frac{1}{-\Delta} (\mathfrak{g})}$. Then, by the identities (\ref{New-force}), for $i=1,2,3$ we write 
\begin{equation}\label{Estim-from-above-2}
\begin{split}
\left| \int_{\Rt} \theta {\bf g}_i  dy \right|=  &\,   \left| \int_{\Rt} \big( \theta {\bf g}_i   - \theta_0{\bf g}_i + \theta_0{\bf g}_i \big)dy  \right| 
\geq \left| \int_{\Rt} \theta_0 {\bf g}_i dy  \right| - \left| \int_{\Rt} (\theta-\theta_0) {\bf g}_i dy \right|\\
\geq &\,  \kappa^2  \left| \int_{\Rt} \tilde{\theta}_0 \mathfrak{ g}_i dy  \right| - \kappa \left| \int_{\Rt} (\theta-\theta_0) \mathfrak{g}_i dy \right|.
\end{split}
\end{equation}
We must estimate the last term on the right-hand side.  To do  this, we recall that $\theta-\theta_0=B_2(\theta, \vu)$; and by the H\"older inequalities we write
\begin{equation*}
\left| \int_{\Rt} (\theta-\theta_0)  \mathfrak{g}_i  dy \right|  \leq  \| \theta - \theta_0 \|_{L^4} \|  \vec{{\mathfrak{ g}}} \|_{L^{4/3}} \leq   \|  B_2(\theta,\vu) \|_{L^4}  \|  \vec{{\mathfrak{ g}}}  \|_{L^{4/3}}. 
\end{equation*} 
Here, we still need to estimate the term $\|  B_2(\theta,\vu) \|_{L^4}$. By the estimate (\ref{Estim-B2}) (with $\kappa\delta_2$ instead of $\delta$; and where $\delta_2$ is  given  in (\ref{delta2}) with the functions $\vec{\mathfrak{f}}$ and $\mathfrak{g}$) we  have 
\begin{equation*}
\|  B_2(\theta,\vu) \|_{L^4}   \leq 3 \kappa \delta_2 C_1(4) \| m_2 \|_{L^{3/2,\infty}} \| \theta \|_{L^4}=  C \kappa \| \theta \|_{L^4}. 
\end{equation*}
Finally, we must estimate the term $\| \theta \|_{L^4}$.  We use the uniform control (\ref{Uniform-control-un}), where we recall that $u_0$ is defined in (\ref{u-u0});  and for $\mathfrak{u}_0=\Big(  (-\Delta)^{-1} \mathfrak{\vf}, (-\Delta)^{-1} \mathfrak{g} \Big)$ we have $u_0= \kappa \mathfrak{u}_0$. Moreover, for the sake of simplicity we shall denote $C'= \left( \sum_{k=0}^{+\infty} \frac{1}{2^k} \right)$. We thus have $\| \theta \|_{L^4} \leq C' \| u_0 \|_{L^4} =C' \kappa \| \mathfrak{u}_0 \|_{L^4}$.  From these estimates we finally obtain 
\begin{equation*}
\left| \int_{\Rt} (\theta-\theta_0)  \mathfrak{g}_i  dy \right|  \leq CC' \kappa^2 \| \mathfrak{u}_0 \|_{L^4}  \|  \vec{{\mathfrak{ g}}}  \|_{L^{4/3}}.
\end{equation*}

We thus get back to the estimate (\ref{Estim-from-above-2})  to write 
\begin{equation*}
\left| \int_{\Rt} \theta {\bf g}_i  dy \right| \geq \kappa^2 \left| \int_{\Rt} \tilde{\theta}_0 \mathfrak{ g}_i dy  \right| - \kappa^3 CC'  \| \mathfrak{u}_0 \|_{L^4}  \|  \vec{{\mathfrak{ g}}}  \|_{L^{4/3}},
\end{equation*}
hence we have 
\begin{equation*}
\kappa^2 \left| \int_{\Rt} \tilde{\theta}_0 \mathfrak{ g}_i dy  \right| - \kappa^3 CC'  \| \mathfrak{u}_0 \|_{L^4}  \|  \vec{{\mathfrak{ g}}}  \|_{L^{4/3}} >0,
\end{equation*}
as long as 
\begin{equation}\label{kappa}
\kappa < \frac{\left| \int_{\Rt} \tilde{\theta}_0 \mathfrak{ g}_i dy  \right|}{CC'  \| \mathfrak{u}_0 \|_{L^4}  \|  \vec{{\mathfrak{ g}}}  \|_{L^{4/3}}} \quad \mbox{and} \quad  0<\left| \int_{\Rt} \tilde{\theta}_0 \mathfrak{ g}_i dy  \right|, 
\end{equation}
where this last inequality is ensured by the assumption (\ref{Cond-g}).  \finpv  

\medskip

Once we have   $M_1>0$,  there exists quantity  $M_2>0$ large enough  such that for $|x|>M_2$ we have 
\[ \left| \textit{o}\left(\frac{1}{|x|}\right)  \right| \leq \frac{M_1}{2|x|}.\]
Then, we get back to the  inequality (\ref{Estim-below}) to obtain the wished estimate from below (\ref{Estim-from-above}).    Proposition \ref{Prop-Estimate-Below} is proven. \finpv 

\medskip

As noticed, Proposition \ref{Prop-Nonexistence} is now proven by Propositions \ref{Prop-Pointwise}, \ref{Prop-Profile} and \ref{Prop-Estimate-Below}. \finpv

\section{Proof of Theorem \ref{Th2}}\label{Sec-Regularity} 
For the sake of clearness, we shall divide the proof  in three main steps.

\medskip

{\bf Step 1. The parabolic system.} Our starting point is the study the time-dependent Boussinesq system. For the velocity $\vv: [0,+\infty) \times \Rt \to \Rt$ and the temperature $\vartheta: [0,+\infty) \times \Rt \to \R$, and moreover, for the time-independent functions $\vf, g, \vec{\bf {g}}$, we consider: 
 \begin{equation}\label{Boussinesq-Parabolic}
\begin{cases}\vspace{2mm}
\partial_t \vv -\Delta \vv + \text{div}(\vv \otimes \vv) + \vec{\nabla} P = \vartheta \vec{{\bf g}}+ \vf, \qquad \text{div}(\vv)=0, \\ \vspace{2mm} 
\partial_t \vartheta-\Delta \vartheta + \text{div}(\vartheta \vv) = g, \\
\vv(0,\cdot) =\vv_0, \ \ \vartheta(0, \cdot)=\vartheta_0,
\end{cases} 
\end{equation}
where $\vv_0: \Rt \to \Rt$, $\vartheta_0:\Rt \to \R$ are the initial data.  For a time $0<T<+\infty$, we denote  $\mathcal{C}_{*}([0,T], L^{p,\infty}(\Rt))$ the functional space of bounded and  weak$-*$ continuous  functions from $[0,T]$ with values in the Lorentz space $L^{p,\infty}(\Rt)$. We prove now the following: 
\begin{Proposition}\label{Prop1} Consider the initial value problem (\ref{Boussinesq-Parabolic}) where $\vf, g$ verify (\ref{Condition-Reg-force}) and $\vec{{\bf g}}$ verifies (\ref{Cond-Reg-g}). Let $p>3$ and let   $\vv_0 \in  L^{p,\infty}(\Rt)$, $ \vartheta_0 \in L^{p,\infty}(\Rt)$  be the initial data. There exists a time $T_0>0$,  depending on $\vv_0$, $\vartheta_0$, $\vf, g$ and $\vec{\bf g}$; and there exists  a unique solution  $(\vv, \vartheta) \in \mathcal{C}_{*}([0,T_0], L^{p,\infty}(\Rt))$ to the equation (\ref{Boussinesq-Parabolic}).  Moreover this solution verifies:   
	\begin{equation}\label{Cond-Sup}
	\sup_{0< t \leq T_0} t^{\frac{3}{2 p}} \Vert  \vv(t, \cdot) \Vert_{L^{\infty}}+	\sup_{0<t \leq T_0} t^{\frac{3}{2 p}} \Vert  \vartheta(t, \cdot) \Vert_{L^{\infty}}  <+\infty. 
	\end{equation}
\end{Proposition} 	
\pv  Mild  solutions of the system (\ref{Boussinesq-Parabolic}) write down as the integral formulation:
\begin{equation}\label{equ}
\begin{split}
\vv(t,\cdot)= & \,\, e^{t \Delta}\vv_0  + \int_{0}^{t} e^{(t-s)\Delta} \,  \P( \vf )ds +\underbrace{ \int_{0}^{t} e^{(t-s)\Delta }\P (div (\vv \otimes \vv ))(s,\cdot) ds}_{\mathcal{B}_1(\vv,\vv)} + \underbrace{\int_{0}^{t} e^{(t-s)\Delta} \P( \theta \vec{\bf g }(s,\cdot)) ds}_{\mathcal{L}(\theta)}, 
\end{split}
\end{equation}
and 
\begin{equation}\label{eqtheta}
\begin{split}
\vartheta(t,\cdot)= & \,\, e^{t\Delta} \vartheta_0 + \int_{0}^{t}e^{(t-s)\Delta} \, g  ds  + \underbrace{ \int_{0}^{t} e^{(t-s)\Delta} \text{div}(\vartheta \vv)(s,\cdot) ds}_{\mathcal{B}_2(\vartheta,\vv)}. 
\end{split}
\end{equation}
By the  Picard's fixed point argument,  we will solve both problems  (\ref{equ}) and (\ref{eqtheta}) in the Banach space  $\mathcal{X}$, where
\[\mathcal{X} = \left\{  f \in \mathcal{C}_{*}([0,T], L^{p,\infty}(\Rt)): \sup_{0<t\leq T} t^{\frac{3}{2 p}} \Vert f(t,\cdot)\Vert_{L^{\infty}}<+\infty \right\}, \]
with the norm 
\[ \Vert f \Vert_{\mathcal{X}}= \sup_{0\leq t \leq T} \Vert f(t,\cdot)\Vert_{L^{p,\infty}}+\sup_{0<t \leq T} t^{\frac{3}{2p}} \Vert f(t,\cdot)\Vert_{L^{\infty}}.\]
Let us mention that for $f_1, f_2 \in \mathcal{X}$, for simplicity we shall write $\| (f_1, f_2) \|_{\mathcal{X}}=\| f_1 \|_{\mathcal{X}}+ \| f_2 \|_{\mathcal{X}}$. 

\medskip 

We start by studying the  terms involving the data in the equations (\ref{equ}) and (\ref{eqtheta}). 
\begin{Lemme}\label{Lem-Data-Linfty}  The next estimate holds:
\begin{equation}\label{Estimate-Data}
\left\Vert \left(e^{t\Delta} \vv_0, e^{t\Delta} \vartheta_0\right) \right\Vert_{\mathcal{X}} + \left\Vert \int_{0}^{t} e^{(t-s)\Delta}\left( \vf,g\right) ds \right\Vert_{\mathcal{X}}  
 \leq  C \Vert (\vv_0, \vartheta_0) \Vert_{L^{p,\infty}} + C T  \left\Vert \left(\vf, g  \right) \right\Vert_{L^{p,\infty}}. 
\end{equation}	
\end{Lemme}	
\pv 
First, for the initial data  $(\vv_0, \varphi)  \in L^{p,\infty}(\Rt)$, by the second point of  Lemma \ref{Lem-Young}  we have the estimate 
$\ds{\Vert (e^{t\Delta} \vv_0, e^{t\Delta} \vartheta_0) \Vert_{L^{p,\infty}} \leq C_2 \Vert (\vv_0, \vartheta_0) \Vert_{L^{p,\infty}}}$.  We thus get $\ds{(e^{t \Delta}\vv_0, e^{t \Delta}\vartheta_0 )  \in \mathcal{C}_{*}([0,T], L^{p,\infty}(\Rt))}$. On the other hand,  by Lemma \ref{Lem-Lorentz-Besov}  we write $\ds{\sup_{0<t\leq T}t^{\frac{3}{2p}} \left\Vert e^{t\Delta} \vv_0\right\Vert_{L^{\infty}}  \leq C \Vert \vv_0 \Vert_{L^{p,\infty}}}$ and $\ds{\sup_{0<t\leq T}t^{\frac{3}{2p}} \left\Vert e^{t\Delta} \vartheta_0\right\Vert_{L^{\infty}}  \leq C \Vert \vartheta_0 \Vert_{L^{p,\infty}}}$. We thus have  $(e^{t \Delta}\vv_0,  e^{t\Delta} \vartheta_0) \in \mathcal{X}$  and the following estimate holds:
\begin{equation}\label{Lin-Initial-Data}
\left\Vert \left(e^{t\Delta} \vv_0, e^{t\Delta} \vartheta_0\right) \right\Vert_{\mathcal{X}} \leq C \Vert (\vv_0, \vartheta_0) \Vert_{L^{p,\infty}}.  
\end{equation} 
 
Thereafter, recall that the external forces $\vf, g$  are time-independent functions. Since we assume (\ref{Condition-Reg-force})   by the  second point of Lemme \ref{Lem-Young} we can write: 
\begin{equation*}
\left\Vert \int_{0}^{t} e^{(t-s)\Delta} \left( \P(\vf),g\right)ds  \right\Vert_{L^{p,\infty}}  \leq  \int_{0}^{t} \left\Vert e^{(t-s)\Delta} \left( \P(\vf), g \right) \right\Vert_{L^{p,\infty}}\, ds 
\leq  C_2 \left\Vert \left( \P(\vf), g \right) \right\Vert_{L^{p,\infty}} \left(\int_{0}^{t} ds\right),
\end{equation*} and we   get
\begin{equation*}\label{estim-forces-1}
\sup_{0 \leq t \leq T} \left\Vert \int_{0}^{t} e^{(t-s)\Delta} \left( \P(\vf),g \right) ds  \right\Vert_{L^{p,\infty}}  \leq C T \left\Vert \left(\P(\vf), g \right) \right\Vert_{L^{p,\infty}}.
\end{equation*} 
On the other hand, to estimate the expression $\ds{\sup_{0<t \leq T} t^{\frac{3}{2p}} \left\| \int_{0}^{t} e^{(t-s)\Delta} \P(\vf) ds \right\|_{L^\infty}}$ we remark that by  Lemma \ref{Lem-Lorentz-Besov} we have
\[ \left\Vert e^{(t-s)\Delta} \P(\vf) \right\Vert_{L^{\infty}} \leq C (t-s)^{-\frac{3}{2p}} \left\Vert \P(\vf)\right\Vert_{L^{p,\infty}},\]
and then we can  write  
\begin{equation*}
\begin{split}
&\, t^{\frac{3}{2 p}}\, \left\Vert \int_{0}^{t} e^{(t-s)\Delta} \P(\vf) ds \right\Vert_{L^\infty}  \leq t^{\frac{3}{2 p}} \, \int_{0}^{t} \left\Vert e^{(t-s)\Delta} \P(\vf) \right\Vert_{L^{\infty}} ds \\
 \leq &\, C\, t^{\frac{3}{2 p}} \, \int_{0}^{t} (t-s)^{-\frac{3}{2p}} \left\Vert \P(\vf) \right\Vert_{L^{p,\infty}} ds  \leq C  t^{\frac{3}{2 p}}  \,  \left\Vert \P(\vf) \right\Vert_{L^{p,\infty}}  \, \left( \int_{0}^{t} (t-s)^{-\frac{3}{2p}} \, ds  \right) \leq C \, t\, \left\Vert \P(\vf)\right\Vert_{L^{p,\infty}}. 
\end{split}
\end{equation*}
We thus obtain
\begin{equation*}
\sup_{0 <t \leq T} t^{\frac{3}{2 p}} \, \left\Vert \int_{0}^{t} e^{(t-s)\Delta} \P(\vf)  ds \right\Vert_{L^\infty}  \leq C T\, \left\Vert  \P(\vf) \right\Vert_{L^{p,\infty}}. 
\end{equation*}
The other term involving the external force $g$ follows similar estimates  and we have 
\begin{equation*}
\sup_{0 <t \leq  T} t^{\frac{3}{2 p}} \, \left\Vert \int_{0}^{t} e^{(t-s)\Delta} g   ds \right\Vert_{L^\infty}  \leq C T \,  \left\Vert g \right\Vert_{L^{p,\infty}}.  
\end{equation*}

By the estimates above  we get
\begin{equation}\label{estim-forces}
\left\Vert \int_{0}^{t} e^{(t-s)\Delta}\left(\P(\vf), g \right) ds \right\Vert_{\mathcal{X}} \leq C T \, \left\Vert \left(\P(\vf), g \right) \right\Vert_{L^{p,\infty}}. 
\end{equation}
Thus, the wished estimate (\ref{Estimate-Data}) follows from the estimates  (\ref{Lin-Initial-Data}) and (\ref{estim-forces}). \finpv 

\medskip

We study now the bilinear terms in (\ref{equ}) and (\ref{eqtheta}). 
\begin{Lemme}\label{Lem-BiLin-Linfty} The following estimates hold:
\begin{equation}\label{Estimate-BiLin}
\| \mathcal{B}_1(\vv, \vv)\|_{\mathcal{X}} \leq C\, T^{\frac{1}{2}-\frac{3}{2p}} \| \vv \|^{2}_{\mathcal{X}}, \quad \| \mathcal{B}_2(\vartheta,\vv)\|_{\mathcal{X}} \leq C\, T^{\frac{1}{2}-\frac{3}{2 p}} \| \vartheta \|_{\mathcal{X}} \|\vv \|_{\mathcal{X}}\, ,
\end{equation}	
where $\frac{1}{2}- \frac{3}{2p}>0$ as long as $p>3$. 
\end{Lemme}	
\pv  First,  the term  $\mathcal{B}_{1}(\vu,\vu)$  in (\ref{equ}) is  estimated as follows: 
\begin{equation}\label{Estim-non-lin-1}
\sup_{0\leq t \leq T} \left\Vert \mathcal{B}_1(\vv, \vv)  \right\Vert_{L^{p,\infty}}  \leq  C\, T^{\frac{1}{2}-\frac{3}{2p}}\, \Vert \vv \Vert^{2}_{E_T}.
\end{equation}
Indeed, by the second point of Lemma \ref{Lem-Young}, by the  well-known estimate on the heat kernel: $\Vert \vec{\nabla}h_{(t-s)} (\cdot)\Vert_{L^1} \leq \frac{c}{(t-s)^{1/2}}$, and moreover, by Lemma \ref{Lem-Lorentz-Besov}  we have
\begin{equation*}\label{BiLin1}
\begin{split}
 \sup_{0\leq t \leq T} \left\Vert \mathcal{B}_1(\vv, \vv)  \right\Vert_{L^{p,\infty}} \leq &\,\,  C\,  \sup_{0\leq t \leq T} \int_{0}^{t}  \left\Vert   e^{(t-s)\Delta } \P \text{div} (\vv \otimes \vv)(s,\cdot)\right\Vert_{L^{p,\infty}}  ds\\
\leq & \,\,C\,  \sup_{0\leq t \leq T} \int_{0}^{t} \frac{1}{(t-s)^{1/2}}  \Vert \vv(s,\cdot) \otimes  \vv(s,\cdot)\Vert_{L^{p,\infty}}  ds\\
\leq &\,\, C\, \sup_{0\leq t \leq T} \int_{0}^{t} \frac{1}{(t-s)^{\frac{1}{2}}\,s^{\frac{3}{2 p}} }  (s^{\frac{3}{2p}}\Vert \vv(s,\cdot) \Vert_{L^{\infty}}) \Vert \vv(s,\cdot)\Vert_{L^{p,\infty}} ds\\
\leq & \,\,C\, T^{\frac{1}{2}-\frac{3}{2p}}\,  \left\Vert  \vv  \right\Vert^{2}_{\mathcal{X}}.  
\end{split}
\end{equation*}  

Thereafter, we have also the estimate 
\begin{equation}\label{Estim-non-lin-2}
\sup_{0 < t \leq T} t^{\frac{3}{2p}} \left\Vert \mathcal{B}_1(\vv, \vv)  \right\Vert_{L^{\infty}} \leq  C\, T^{\frac{1}{2}-\frac{3}{2p}}\,  \left\Vert \vv  \right\Vert^{2}_{\mathcal{X}}. 
\end{equation}
Indeed, we write:
\begin{equation*}
\begin{split}
\sup_{0 < t \leq T} t^{\frac{3}{2p}} \left\Vert \mathcal{B}_1(\vv, \vv)  \right\Vert_{L^{\infty}}  \leq  \,\,  \sup_{0<t<T} t^{\frac{3}{2p}}    \int_{0}^{t}  \left\Vert  e^{(t-s)\Delta } \P  \text{div} \left( \vv \otimes \vv  \right)(s,\cdot)  \right\Vert_{L^{\infty}} ds=(a).
\end{split}
\end{equation*}
Here, we recall that the operator $e^{(t-s)\Delta} \P(div(\cdot))$ writes down as  a tensor of convolution operators (in the spatial variable) whose kernels $K_{i,j,k}$ verify $\ds{\vert K_{i,j,k}(t-s, x )\vert \leq \frac{c}{((t-s)^{1/2}+\vert x \vert)^{4}}}$. See \cite[Proposition $11.1$]{PGLR}.   Then, we have  $\ds{\Vert K_{i,j}(t-s,\cdot) \Vert_{L^1} \leq \frac{c}{(t-s)^{1/2}}}$; and we can write:
\begin{equation*}
\begin{split}
(a)\leq &  \,\, C\,  \sup_{0 < t \leq T} t^{\frac{3}{2p}}   \int_{0}^{t} \frac{1}{(t-s)^{1/2}} \left( \Vert \vv(s,\cdot) \otimes  \vv(s,\cdot)\Vert_{L^{\infty}} \right) ds\\
\leq & \,\,  C\,  \sup_{0 < t \leq T} t^{\frac{3}{2p}}   \int_{0}^{t} \frac{ds}{(t-s)^{1/2} s^{\frac{3}{p}}} \left( s^{\frac{3}{2p}} \Vert \vv(s,\cdot) \Vert_{L^{\infty}}\right)^2 ds\\
\leq & \,\,  C\, \left( \sup_{0 < t \leq T} t^{\frac{3}{2p}}    \int_{0}^{t} \frac{ds}{(t-s)^{1/2} s^{\frac{3}{p}}} \right) \left\Vert  \vv \right\Vert^{2}_{\mathcal{X}}. \\
\leq &  \,\,  C\, T^{\frac{1}{2}-\frac{3}{2p}}\,  \left\Vert \vv  \right\Vert^{2}_{\mathcal{X}}.
\end{split}
\end{equation*}
By the estimates (\ref{Estim-non-lin-1}) and (\ref{Estim-non-lin-2}) we obtain
\begin{equation}\label{Estim_B1}
\| \mathcal{B}_1 (\vv, \vv) \|_{\mathcal{X}} \leq c\, T^{\frac{1}{2}-\frac{3}{2 p}}\, \| \vv \|^{2}_{\mathcal{X}}.
\end{equation}

We estimate now the term $\mathcal{B}_2(\vartheta,\vv)$. By following the same arguments in the prove of the estimate (\ref{Estim-non-lin-1}) (where we use the information $s^{\frac{3}{2 p}} \| \vv(s,\cdot)\|_{L^\infty} <+\infty$  and $\| \vartheta(s,\cdot)\|_{L^{p,\infty}} <+\infty$ ) we have 
\begin{equation}\label{Estim-non-lin-3}
\sup_{0 \leq t \leq T} \| \mathcal{B}_2(\vartheta,\vv) \|_{L^{p,\infty}} \leq C\, T^{\frac{1}{2}-\frac{3}{2p}} \,  \| \vartheta \|_{\mathcal{X}} \| \vv \|_{\mathcal{X}}. 
\end{equation}
Moreover, by following  similar  ideas  in the proof of the estimate (\ref{Estim-non-lin-2}) we can prove: 
\begin{equation}\label{Estim-non-lin-4}
\sup_{0<t \leq T} t^{\frac{3}{2 p}}\| \mathcal{B}_2(\vartheta,\vv)\|_{L^\infty} \leq C\, T^{\frac{1}{2}-\frac{3}{2 p}}  \| \vartheta \|_{\mathcal{X}} \| \vv \|_{\mathcal{X}} .
\end{equation}
Indeed,  we just  write 
\begin{equation*}
\sup_{0<t \leq T} t^{\frac{3}{2 p}}\| \mathcal{B}_2(\vartheta,\vv)\|_{L^\infty} \leq \, C \, \left( \sup_{0<t\leq T} t^{\frac{3}{2 p}} \int_{0}^{t} \frac{ds}{(t-s)^{1/2} s^{\frac{3}{2p}+\frac{3}{2p}}}\right)\| \vartheta \|_{\mathcal{X}} \| \vv \|_{\mathcal{X}}\,   \leq\, C\, T^{\frac{1}{2}-\frac{3}{2 p}}\, \| \vartheta \|_{\mathcal{X}} \| \vv \|_{\mathcal{X}}.
\end{equation*} 
By the estimates (\ref{Estim-non-lin-3}) and (\ref{Estim-non-lin-4}) we get 
\begin{equation}\label{Estimate_B2}
\| \mathcal{B}_2(\vartheta,\vv) \|_{\mathcal{X}} \leq C\, T^{\frac{1}{2}-\frac{3}{2p}}\, \| \vartheta \|_{\mathcal{X}} \| \vv \|_{\mathcal{X}}.
\end{equation}

To finish the proof, we gather the estimates (\ref{Estim_B1}) and (\ref{Estimate_B2}) to obtain the wished estimate (\ref{Estimate-BiLin}). \finpv

\medskip 

Finally, we must study the linear term $\mathcal{L}(\vartheta)$ in the equation (\ref{equ}). Recall that by the assumption (\ref{Cond-Reg-g}) we have $\vec{\bf g } \in L^\infty(\Rt)$. 
\begin{Lemme}\label{Lem-Lin-Linfty}  The following estimate holds:
	\begin{equation}\label{Estimate-Lin}
	\| \mathcal{L}(\vartheta) \|_{\mathcal{X}} \leq C \| \vec{\bf g } \|_{L^\infty} T \, \|  \vartheta \|_{\mathcal{X}}. 
	\end{equation}
\end{Lemme}	
\pv We start with the following estimate, where by the second point of  Lemma \ref{Lem-Young} and the H\"older inequalities we write
\begin{equation*}
\begin{split}
&\sup_{0\leq t \leq T} \, \left\|  \int_{0}^{t} e^{(t-s)\Delta} \P (\vartheta \vec{\bf g})(s,\cdot) ds\right\|_{L^{p,\infty}} \leq  \sup_{0\leq t \leq T} \, \int_{0}^{t} \left\| e^{(t-s)\Delta}  \P (\vartheta \vec{\bf g})(s,\cdot) \right\|_{L^{p,\infty}} ds \\
\leq &\,  C\,  \sup_{0\leq t \leq T} \,  \int_{0}^{t} \| \vartheta \vec{\bf g }(s,\cdot)\|_{L^{p,\infty}} ds  \leq  C  \| \vec{\bf g}   \|_{L^\infty} T \| \vartheta \|_{\mathcal{X}}. 
\end{split}
\end{equation*}
Moreover,  by the third point of Lemma \ref{Lem-Young}  and by using again the H\"older inequalities,   the next estimate holds
\begin{equation*}
\begin{split}
& \sup_{0 < t \leq  T} \, t^{\frac{3}{2p}} \,  \left\|  \int_{0}^{t} e^{(t-s)\Delta} \P (\theta \vec{\bf g})(s,\cdot) ds\right\|_{L^{\infty}}  \leq  C\, \sup_{0 < t \leq  T} \, t^{\frac{3}{2p}} \, \int_{0}^{t} (t-s)^{-\frac{3}{2 p}} \| \P(\vartheta \vec{\bf g })(s,\cdot)\|_{L^{p,\infty}} ds \\
\leq &\, C\, \left(  \sup_{0 < t \leq  T} \, t^{\frac{3}{2p}} \, \int_{0}^{t} (t-s)^{-\frac{3}{2 p}} ds \right) \, \| \vartheta \vec{\bf g } \|_{\mathcal{X}} \leq C  \| \vec{\bf g} \|_{L^\infty}  T\|  \vartheta \|_{\mathcal{X}}. 
\end{split}
\end{equation*}
We thus  obtain the wished estimate (\ref{Estimate-Lin}). \finpv

\medskip

With the Lemmas \ref{Lem-Data-Linfty}, \ref{Lem-BiLin-Linfty} and \ref{Lem-Lin-Linfty} at our disposal, we set a time $T_0>0$  small enough and Proposition \ref{Prop1} follows from standard arguments.  \finpv

\medskip

{\bf Step 2. Global boundness of $\vu$ and $\theta$}. Now, we get back to steady-state system (\ref{Boussinesq}). With the help of Proposition \ref{Prop1},  we are able to prove  the following:  
\begin{Proposition}\label{Prop2}    Let  $(\vu,P,\theta)$ be a  weak solution of the system  (\ref{Boussinesq}) given in  Definition \ref{Def-weakSol-Boussinesq}. If $\vu \in L^{p,\infty}(\Rt)$ and $\theta \in  L^{p,\infty}(\Rt)$   with $p>3$, then we have $\vu \in L^{\infty}(\Rt)$ and $\theta\in L^{\infty}(\Rt)$.  
\end{Proposition} 	
\pv   In the initial value  problem (\ref{Boussinesq-Parabolic}), we set  the initial data $\ds{(\vv_0, \vartheta_0)=(\vu, \theta)}$. Then, by  Proposition \ref{Prop1} there exists a time $0<T_0$ and there exists  a unique arising  solution $(\vv, \vartheta) \in \mathcal{C}_{*}([0,T_0], L^{p,\infty}(\Rt))$ to the equation (\ref{Boussinesq-Parabolic}). 

\medskip 

On the other hand, we have the following key remark. Since   $\vu$ and $ \theta$ are time-independent functions we have $\partial_t \vu=0$ and $\partial_t \theta=0$. Thus, the  couple $\ds{(\vu,\theta )}$ is also a solution of the initial value problem (\ref{Boussinesq-Parabolic}) with the same  initial data $(\vu, \theta)$, and moreover  we have $(\vu, \theta) \in \mathcal{C}_{*}([0,T_0], L^{p,\infty}(\Rt))$. 
 
\medskip

Consequently, in the space $\mathcal{C}_{*}([0,T_0], L^{p,\infty}(\Rt))$ we  have   two solutions of  (\ref{Boussinesq-Parabolic}) with the same initial data: on the one hand, the solution $\ds{(\vv , \vartheta )}$ given by  Proposition \ref{Prop1} and, on the other hand,  the (steady-sate) solution $\ds{(\vu, \theta)}$. By  uniqueness we have the identity $\ds{(\vv, \vartheta)= (\vu, \theta)}$ and by   (\ref{Cond-Sup})  we can write 
\[  \sup_{0<t \leq T_0} t^{\frac{3}{2 p}} \Vert  \vu \Vert_{L^{\infty}}+\sup_{0<t \leq T_0} t^{\frac{3}{2 p}}  \Vert  \theta \Vert_{L^{\infty}} <+\infty. \]
But, as the solution $ (\vu , \theta)$ does not depend on the time variable we  have $\vu  \in L^{\infty}(\Rt)$ and $\theta  \in L^{\infty}(\Rt)$.  Proposition  \ref{Prop2} is now proven. \finpv     

\medskip

{\bf Step 3. Estimates on high order derivatives in  Lorentz spaces}. The global boundness of $\vu$ and $\theta$ yields the following: 
\begin{Proposition}\label{Prop3} For $k \geq 0$ assume that $\vf, g$ verify (\ref{Condition-Reg-force}) and $\vec{\bf g}$ verifies (\ref{Cond-Reg-g}).  Moreover, assume that $\vu \in L^{p,\infty}(\Rt)$ and  $\theta \in  L^{p,\infty}(\Rt)$ with  $p>3$.  Then  we have $\vu \in \mathcal{ W}^{k+2,p}(\Rt)$, $\theta\in \mathcal{W}^{k+2,p}(\Rt)$ and  $P \in \mathcal{W}^{k+1,p}(\Rt)$.    
\end{Proposition}	 
\pv We will use the  integral formulations (\ref{Boussineq-Int-1}) and (\ref{Boussineq-Int-2}) to prove  that for all multi-indice $\vert \alpha \vert \leq k+2$  and for all  $1\leq \sigma <+\infty$ we have  $\partial^{\alpha} \vu  \in L^{p \sigma,\infty}(\Rt)$  and $\partial^{\alpha} \theta \in L^{p \sigma,\infty}(\Rt)$. We shall prove this fact by an  iteration process  respect to the order of the multi-indices $\alpha$, which  we will denote as $\vert \alpha \vert$. For the  sake of clearness,  in the following couple of  technical lemmas  we separately  prove each step in the iterative argument.
\begin{Lemme}[The case initial case]  Recall that by Proposition \ref{Prop2} we have   $\vu , \theta  \in L^{\infty}(\Rt)$. Then,  for $\vert \alpha \vert \leq 2$ and for all $1\leq \sigma < +\infty$ we have $\partial^{\alpha}\vu  \in L^{p \sigma,\infty}(\Rt)$ and $\partial^{\alpha}\theta \in L^{q \sigma,\infty}(\Rt)$. 
\end{Lemme}	 
\pv Due to the  coupled structure of the  equations (\ref{Boussineq-Int-1})-(\ref{Boussineq-Int-2}),   we  must study first the equation (\ref{Boussineq-Int-2}) and  then we  study  the equation (\ref{Boussineq-Int-1}). 

\begin{itemize}
\item Let $|\alpha|=1$.  By the equation (\ref{Boussineq-Int-2}) we write 
\begin{equation}\label{Der-theta}
\partial^{\alpha} \theta = -\frac{\partial^\alpha}{-\Delta}(\text{div} (\theta \vu)) + \frac{\partial^\alpha}{-\Delta}(g). 
\end{equation}
Then,  remark  that since $\theta \in L^{p,\infty}(\Rt) \cap L^\infty(\Rt)$ by the first point of Lemma \ref{Lem-Interpolation} we have $\theta \in L^{p \sigma,\infty}(\Rt)$ for all $1\leq \sigma < \infty$. Moreover,  since $\vu \in L^\infty(\Rt)$ we have $\theta \vu \in L^{p\sigma,\infty}(\Rt)$.  With this information, for the first term on the right-hand side, remark that  since $|\alpha|=1$ the operator $\ds{-\frac{\partial^\alpha}{-\Delta}(\text{div} (\cdot))}$ writes down as a combination of the Riesz Transforms $\mathcal{R}_{i}\mathcal{R}_{j}$ (with $i,j=1,2,3$).  Thus, as we have $\theta \vu \in L^{p \sigma,\infty}(\Rt)$, by Lemma \ref{Lem-Riesz} we obtain 
$\ds{-\frac{\partial^\alpha}{-\Delta}(\text{div} (\theta \vu))} \in L^{p \sigma, \infty}(\Rt)$. 

\medskip

 On the other hand, for the second term on the right-hand side, recall that the operator $\ds{\frac{\partial^\alpha}{-\Delta}(\cdot)}$ writes down as a product of convolution with the tensor $m_2 \in L^{3/2,\infty}(\Rt)$. Moreover, recall that by the assumption (\ref{Condition-Reg-force})  we have $g \in L^{\frac{3p}{3+p},\infty}(\Rt)$. We thus apply the  first point of Lemma \ref{Lem-Young}  (with $1+1/p= 2/3+ (3+p)/3p$) to obtain  $\ds{\frac{\partial^\alpha}{-\Delta}(g) \in L^{p,\infty}(\Rt)}$. Thereafter, also by  (\ref{Condition-Reg-force}) and by the second point of Lemma \ref{Lem-Interpolation} we have $g\in L^{3,1}(\Rt)$ (remark that $\frac{3p}{3+p} <3$) and by the third point of Lemma \ref{Lem-Young} we get $\ds{\frac{\partial^\alpha}{-\Delta}(g) \in L^{\infty}(\Rt)}$. Finally, by the first point of Lemma \ref{Lem-Interpolation} we have $\ds{\frac{\partial^\alpha}{-\Delta}(g) \in L^{p\sigma,\infty}(\Rt)}$.  Consequently, we get $\ds{\partial^{\alpha} \theta} \in L^{p \sigma, \infty}(\Rt)$ for all $1\leq \sigma < +\infty$. 

\medskip

With this information, we can study now the equation (\ref{Boussineq-Int-1}). As before,   we write 
\begin{equation}\label{Der-u}
\partial^\alpha \vu =  - \frac{\partial^\alpha}{-\Delta} \left( \P ( \text{div} (\vu \otimes \vu) ) \right)+ \frac{\partial^\alpha}{-\Delta}\Big( \P \left(  \theta \vec{\bf g} \right)\Big)+  \frac{\partial^\alpha}{-\Delta}\Big( \P (\vf)\Big),
\end{equation} 
and we shall prove that each term on right-hand side belongs to the space $L^{p \sigma, \infty}(\Rt)$ for all $1\leq \sigma < +\infty$. Indeed, for the first term we recall that since $\vu \in L^{p,\infty}(\Rt)\cap L^{\infty}(\Rt)$ by Lemma \ref{Lem-Interpolation} we have $\vu \in L^{p \sigma, \infty}(\Rt)$;  and consequently $\vu \otimes \vu \in L^{p \sigma, \infty}(\Rt)$. Then, we follow the same ideas above to conclude that $\ds{- \frac{\partial^\alpha}{-\Delta} \left( \P ( \text{div} (\vu \otimes \vu) ) \right) \in L^{p \sigma, \infty}(\Rt)}$. 

\medskip

For the second term, recall that the operator $\ds{\frac{\partial^\alpha}{-\Delta} \P(\cdot)}$ writes down as a product of convolution with the tensor $m_1 \in L^{3/2,\infty}(\Rt)$. Then, by the first point of Lemma \ref{Lem-Young} we write $\left\|  \frac{\partial^\alpha}{-\Delta} \P(\theta \vec{\bf g}) \right\|_{L^{p,\infty}} \leq C_1 \| m_1 \|_{L^{3/2,\infty}} \| \theta \vec{\bf g} \|_{L^{\frac{3p}{3+p}, \infty}}$. Then, by the assumption (\ref{Cond-Reg-g}) and the second point of Lemma \ref{Lem-Holder}  we get $\| m_1 \|_{L^{3/2,\infty}} \| \theta \vec{\bf g} \|_{L^{\frac{3p}{3+p}, \infty}} \leq C \| \theta \|_{L^{p,\infty}} \| \vec{\bf g} \|_{L^{3,\infty}}$. We thus get $\ds{\frac{\partial^\alpha}{-\Delta} \P(\theta \vec{\bf g}) \in L^{p,\infty}(\Rt)}$. On the other hand, by the third point of Lemma \ref{Lem-Young} we can write $\| \frac{\partial^\alpha}{-\Delta} \P(\theta \vec{\bf g})  \|_{L^\infty} \leq C_3 \| m_1 \|_{L^{3/2,\infty}} \|\theta \vec{\bf g} \|_{L^{3,1}} \leq \| \theta \|_{L^\infty}\| \vec{\bf g} \|_{L^{3,1}}$. We thus get $\ds{\frac{\partial^\alpha}{-\Delta} \P(\theta \vec{\bf g})  \in L^{\infty}(\Rt)}$; and consequently we have $\ds{\frac{\partial^\alpha}{-\Delta} \P(\theta \vec{\bf g})  \in L^{p\sigma,\infty}(\Rt)}$.

\medskip 

The first term is similarly estimated as the term $\ds{\frac{\partial^\alpha}{-\Delta}(g)}$ and we have $\ds{\frac{\partial^\alpha}{-\Delta}\Big( \P (\vf)\Big) \in L^{p\sigma,\infty}(\Rt)}$. We thus obtain $\ds{\partial^\alpha \vu \in L^{p\sigma,\infty}(\Rt)}$ for all $1 \leq \sigma < +\infty$. 

\item Let $|\alpha |=2$.  We get back to the expression $\partial^\alpha \theta$ defined in (\ref{Der-theta}). For the first term on right-hand side  we write 
\begin{equation*}
-\frac{\partial^{\alpha}}{-\Delta}(\text{div}(\theta \vu))= -\frac{\partial^{\alpha_1}}{-\Delta}(\text{div}(\partial^{\alpha_2} (\theta \vu))), \quad  \mbox{where} \  \ \alpha=\alpha_1+\alpha_2, \  \ |\alpha_1 |=|\alpha_2|=1. 
\end{equation*}
Here, we shall verify that $\partial^{\alpha_2} (\theta \vu) \in L^{p\sigma,\infty}(\Rt)$. Indeed, for $i=1,2,3$ we write $\partial^{\alpha_2}(\theta u_i)= (\partial^{\alpha_2}\theta) u_i + \theta (\partial^{\alpha_2} u_i)$ and since we get $\partial^{\alpha_2}\theta, \partial^{\alpha_2} \vu \in L^{p\sigma, \in}$, $\vu, \theta \in L^{\infty}(\Rt)$ we directly obtain    $\partial^{\alpha_2} (\theta \vu) \in L^{p\sigma,\infty}(\Rt)$. Thereafter, recall  that the operator $\ds{-\frac{\partial^{\alpha_1}}{-\Delta}(\text{div}(\cdot))}$ writes down as a combination of the Riesz transforms $\mathcal{R}_{i}\mathcal{R}_j$  and by applying the  Lemma \ref{Lem-Riesz} we have $\ds{-\frac{\partial^{\alpha}}{-\Delta}(\text{div}(\theta \vu)) \in L^{p\sigma,\infty}(\Rt)}$.  The second term on the right-hand side of the expression  (\ref{Der-theta}) also belongs to the space $L^{p\sigma,\infty}(\Rt)$ due to the hypothesis (\ref{Condition-Reg-force}). \\

We study now the expression $\partial^{\alpha} \vu$ defined in (\ref{Der-u}). The first term on the right-hand side follows the same arguments above and we have $\ds{- \frac{\partial^\alpha}{-\Delta} \left( \P ( \text{div} (\vu \otimes \vu) ) \right) \in L^{p\sigma,\infty}(\Rt)}$.  To study the second term on the right-hand side, we also remark that in this case the operator $\ds{\frac{\partial^\alpha}{-\Delta}(\P(\cdot))}$ writes down as a linear combination of the Riesz transforms, and moreover, since $\vec{\bf g} \in L^\infty(\Rt)$ we have $\ds{\frac{\partial^\alpha}{-\Delta}(\P(\theta \vec{\bf g})) \in L^{p\sigma,\infty}(\Rt)}$. Finally, always by the hypothesis (\ref{Condition-Reg-force}) we have $\ds{ \frac{\partial^\alpha}{-\Delta}\Big( \P (\vf)\Big) \in L^{p\sigma,\infty}(\Rt)}$.  \finpv 
\end{itemize} 

\begin{Lemme}[The iterative process] For all $1\leq m \leq k$ (with $k\geq 2$) and for all multi-indice $|\alpha | \leq k$,  assume that $\partial^\alpha \vu \in L^{p\sigma,\infty}(\Rt)$ and $\partial^\alpha \theta \in L^{p\sigma,\infty}(\Rt)$, for all $1\leq \sigma <+\infty$. Then it holds true for all multi-indice $|\alpha | \leq k+2$. 
\end{Lemme} 
\pv  We shall similar ideas in the proof of the previous lemma. 
\begin{itemize}
\item  Let $|\alpha |=k+1$. We start by studying the expression $\partial^{\alpha}\theta$ defined in (\ref{Der-theta}). For the first term on the right-hand side we write 
\begin{equation*}
-\frac{\partial^\alpha}{-\Delta}(\text{div} (\theta \vu))=- \frac{\partial^{\alpha_1}}{-\Delta} (\text{div} (\partial^{\alpha_2}(\theta \vu))), \quad \alpha=\alpha_1+\alpha_2, \ \ |\alpha_1 |=1, \ \ |\alpha_2|=k.
\end{equation*}	
In this last expression, we must verify that $\partial^{\alpha_2}(\theta \vu) \in L^{p\sigma,\infty}(\Rt)$. Indeed, for $i=1,2,3$ by the Leibniz rule we write $\ds{\partial^{\alpha_2}(\theta u_i) = \sum_{|\beta | \leq k} c_{\alpha_2, \beta} \partial^\beta \theta \, \partial^{\alpha_2-\beta} u_i}$, where the constant $c_{\alpha_2,\beta}>0$ depends on the multi-indices $\alpha_2$ and $\beta$. Recall that by hypothesis we have $\partial^\beta \theta, \, \partial^{\alpha_2-\beta} u_i \in L^{p\sigma, \infty}(\Rt)$ and by the H\"older inequalities we get $\partial^\beta \theta \, \partial^{\alpha_2-\beta} u_i \in L^{p\sigma,\infty}(\Rt)$. We thus have $\partial^{\alpha_2}(\theta \vu) \in L^{p\sigma,\infty}(\Rt)$ and consequently $\ds{-\frac{\partial^\alpha}{-\Delta}(\text{div} (\theta \vu)) \in L^{p\sigma,\infty}(\Rt)}$.  For the second term on the right-hand side, we write
\begin{equation*}
\frac{\partial^{\alpha}}{-\Delta}(g)=\frac{\partial^{\alpha_1}}{-\Delta} (\partial^{\alpha_2} g), \quad \mbox{where} \ \ \alpha=\alpha_1+\alpha_2, \ \  |\alpha_1|=2, \ \ |\alpha_2|=k-1, 
\end{equation*}
and by the assumption (\ref{Condition-Reg-force}) we get $\ds{\frac{\partial^\alpha}{-\Delta}(g) \in L^{p\sigma,\infty}(\Rt)}$.

\medskip

Once we have $\partial^{\alpha} \theta \in L^{p\sigma,\infty}(\Rt)$ for all $|\alpha | \leq k+1$, by the expression (\ref{Der-u}) and by following the same arguments above we obtain $\partial^{\alpha}\vu \in L^{p\sigma,\infty}(\Rt)$ for all $|\alpha | \leq k+1$. We just mention  that to treat the term $\ds{\frac{\partial^\alpha}{-\Delta}\Big( \P \left(  \theta \vec{\bf g} \right)\Big)}$ we split 
\[  \frac{\partial^\alpha}{-\Delta}\Big( \P \left(  \theta \vec{\bf g} \right)\Big) = \frac{\partial^\alpha_1}{-\Delta}\Big( \P \left(  \partial^{\alpha_2}( \theta \vec{\bf g}) \right)\Big),  \ \  |\alpha_1|=2, \ \ |\alpha_2|=k-1, \]
where we use again the Leibniz rule to study  the term $\partial^{\alpha_2}( \theta \vec{\bf g})$. 
\item Let $|\alpha | = k+2$.  With the information $\partial^{\alpha} \theta \in L^{p\sigma,\infty}(\Rt)$ and $\partial^{\alpha}\vu \in L^{p\sigma,\infty}(\Rt)$  (for all $|\alpha | \leq k+1$) at our disposal, we just repeat again the arguments above to get $\partial^{\alpha} \theta \in L^{p\sigma,\infty}(\Rt)$ and $\partial^{\alpha}\vu \in L^{p\sigma,\infty}(\Rt)$ for $|\alpha | = k+2$. \finpv 
\end{itemize}	

\begin{Remarque}\label{Rmk4} By the hypothesis (\ref{Condition-Reg-force}) on the external forces $\vf$ and $g$, and moreover, by the terms  $\frac{\partial^\alpha}{-\Delta}\P(\vf)$ and  $\frac{\partial^\alpha}{-\Delta}(g)$ in the  expressions (\ref{Der-u}) and (\ref{Der-theta}) respectively, we get that $k+2$ is the maximum gain of regularity expected for $\partial^{\alpha}\vu$ and $\partial^\alpha \theta$.  
\end{Remarque}	

It remains to study the pressure term $P$ in the first equation of the coupled system (\ref{Boussinesq}). To do this, we apply the divergence operator to this equations to obtain that $P$ is related to the velocity $\vu$, the temperature $\theta$ and the  force $\vf$  through the expression
\begin{equation}\label{Pressure}
P= \frac{1}{-\Delta} \text{div}(\text{div}(\vu \otimes \vu)) - \frac{1}{-\Delta} (\text{div}(\theta \vec{\bf g })) - \frac{1}{-\Delta} (\text{div}(\vf)).  
\end{equation}
By following the same arguments in the proof the previous lemmas, we can prove the next one:
\begin{Lemme}[The pressure term] For all multi-indice $|\alpha|\leq k+1$ we have $\partial^{\alpha} P \in L^{p\sigma,\infty}(\Rt)$. 
\end{Lemme}	
\pv  Let $|\alpha | \leq k+1$. By the expression (\ref{Pressure}) we write 
\begin{equation}\label{Der-Pressure}
\partial^{\alpha}P= \frac{\partial^{\alpha}}{-\Delta} \text{div}(\text{div}(\vu \otimes \vu)) - \frac{\partial^{\alpha}}{-\Delta} (\text{div}(\theta \vec{\bf g})) - \frac{\partial^{\alpha}}{-\Delta} (\text{div}(\vf)), 
\end{equation}
where  we must study each term on the right-hand side. For the first term, we recall that $\partial^{\alpha} \vu \in L^{p\sigma,\infty}(\Rt)$ for all $1\leq \sigma <+\infty$. Then, we can use  the H\"older inequalities as well as  the Leibniz rule to obtain $\ds{\frac{\partial^{\alpha}}{-\Delta} \text{div}(\text{div}(\vu \otimes \vu))= \frac{1}{-\Delta} \text{div}(\text{div}(\partial^{\alpha}(\vu \otimes \vu)))\in L^{p\sigma,\infty}(\Rt)}$.  The second term is similarly treated  by using now the information $\partial^{\alpha} \theta \in L^{p\sigma,\infty}(\Rt)$ and the assumption (\ref{Cond-Reg-g}). Finally, always by the hypothesis (\ref{Condition-Reg-force}) the third term verifies $\ds{\frac{\partial^{\alpha}}{-\Delta} (\text{div}(\text{div}(\mathbb{F}))) \in L^{p\sigma,\infty}(\Rt)}$.  \finpv

\medskip

Now, we are able to finish the proof of Theorem \ref{Th1}. For this we recall the Lorentz space $L^{p,\infty}(\Rt)$ embeds in the Morrey space $\dot{M}^{1,p}(\Rt)$ defined by the expression (\ref{Morrey}). Consequently, for all multi-indice $|\alpha|\leq k+1$, by Lemma \ref{Lem-Holder-Reg} the functions $\partial^{\alpha}\vu$ and $\partial^{\alpha}\theta$ are H\"older continuous with parameter $s=1-3/p$. Moreover, for $|\alpha | \leq k$ this fact also holds true for the function $\partial^\alpha P$. Theorem \ref{Th2} is now proven. \finpv 

\section{Proof of Theorem \ref{Th3}} 
This result is based on the following Caccioppoli-type estimate for the temperature $\theta$. For $R>0$ we shall denote  $B_R=\{ x \in \Rt: \ | x | <R \}$ and $C(R/2, R)=\{ x \in \Rt: \ R/2 < |x | < R \}$. 
\begin{Proposition}\label{Caccioppoli}  Let $(\vu, P, \theta) \in \mathcal{C}^{2}_{loc}(\Rt) \times \mathcal{C}^{1}_{loc}(\Rt) \times \mathcal{C}^{2}_{loc}(\Rt)$ be a smooth solution of the coupled system (\ref{Boussinesq-homog}). Let $p>3$. There exists a constant $C>0$ such that for all $R \geq 1$ the following estimate hold:
\begin{equation}\label{Caccioppoli-theta}
\int_{B_{R/2}} | \vec{\nabla} \theta |^2 dx \leq  C\,   \left( \int_{C(R/2, R)} | \theta |^p dx \right)^{2/p} \left[  R^{1-\frac{6}{p}}  + R^{2-\frac{9}{p}}\, \left( \int_{C(R/2,R)} | \vu  |^p dx \right)^{1/p}\right].
\end{equation}	
\end{Proposition}	
\pv   To prove the estimate (\ref{Caccioppoli-theta}), we introduce  the following  cut-off function: let $\varphi \in \mathcal{C}^{\infty}_{0}(\Rt)$ be a positive and radial function verifying $\varphi (x)=1$ when $|x|<1/2$ and $\varphi(x)=0$  when $|x| \geq 1$. Then, for $R\geq 1$ we define $\varphi_R(x)=\varphi(x  / R)$. Consequently,  the cut-off function $\varphi_R$ satisfies the next useful properties: $\varphi_R(x)=1$ when $|x|< R/2$, $\varphi_2 (x)=0$ when $|x| \geq R$, and moreover, for all multi-indice $\alpha$ we have $\text{supp}(\partial^{\alpha} \varphi_R) \subset C(R/2, R)$. 

\medskip

Now, we multiply each term in the second equation in (\ref{Boussinesq-homog}) by $\varphi_R \theta$; and we integrate on the ball $B_R$  (remark that $\text{supp}(\varphi_R \theta) \subset B_R$) to get 
\begin{equation}\label{Iden-1}
\int_{B_R} (-\Delta \theta )  \varphi_R \theta dx + \int_{B_R} \text{div}(\theta \vu) \varphi_R \theta dx =0. 
\end{equation}  
Here it is worth highlighting each term above is well-defined since we have $\vu, \theta \in \mathcal{C}^{2}_{loc}(\Rt)$. To study the first term we integrate by parts to write
\begin{equation}\label{Grad-theta}
\begin{split}
&\int_{B_R} (-\Delta \theta )  \varphi_R \theta dx  =    -\sum_{i=1}^{3} \int_{B_R} (\partial^{2}_{j} \theta) ( \varphi_R \theta) dx= \sum_{i,j=1}^{3} \int_{B_R} \partial_{i} \theta  \partial_{i}(\varphi_R \theta) dx \\  
=& \,  \sum_{i=1}^{3} \int_{B_R} (\partial_i \theta) (\partial_i  \varphi_R)   \theta  dx +
\sum_{i=1}^{3} \int_{B_R} (\partial_i \theta)  \varphi_R (  \partial_i \theta) dx =  \sum_{i=1}^{3}   \int_{B_R} (\partial_i \varphi_R) (\partial_i \theta)  \theta dx +  \sum_{i=1}^{3} \int_{B_R} \varphi_R  (\partial_i \theta)^2  dx\\ 
=& \,  \sum_{i=1}^{3} \int_{B_R} (\partial_i \varphi_R )  \partial_{i} \left( \frac{\theta^{2}}{2}\right)dx + \int_{B_R}\varphi_R \vert \vec{\nabla}\theta \vert^{2} dx
= - \int_{B_R} \Delta \varphi_R  \left( \frac{\theta^2}{2}  \right)dx +   \int_{B_R}\varphi_R \vert \vec{\nabla} \theta\vert^{2} dx. 
\end{split}
\end{equation}
Thereafter, to study the second term, we remark first that since $\text{div}(\vu)=0$ we have $\text{div}(\theta \vu)= \vu \cdot \vec{\nabla} \theta$. Then,  we use again the integration by parts to get
\begin{equation*}
\int_{B_R} \vu \cdot \vec{\nabla} \theta\cdot ( \varphi_R \theta )dx =   \sum_{i=1}^{3} \int_{B_R} u_i (\partial_i \theta) ( \varphi_R \theta) dx =  \sum_{i=1}^{3} \int_{B_R} \varphi_R u_i  (\partial_i \theta) \theta dx =   \sum_{i=1}^{3} \int_{B_R} \varphi_R u_i  \,  \partial_i \left(\frac{\theta^{2}}{2}\right) dx,
\end{equation*}
but, always by the fact that $div(\vu)=0$,  we can write 
\begin{equation} \label{Transport-theta}
\sum_{i=1}^{3} \int_{B_R} \varphi_R u_i  \,  \partial_i \left(\frac{\theta^{2}}{2}\right) dx =\sum_{i=1}^{3} \int_{B_R} \varphi_R \partial_{i} \left( u_i \frac{\theta^{2}}{2} \right)dx=  - \int_{B_R} \vec{\nabla} \varphi_R \cdot \left( \frac{ \theta^2}{2} \vu\right)dx.  
\end{equation}
With the identities (\ref{Grad-theta}) and (\ref{Transport-theta}) at our disposal, we get back to the identity  (\ref{Iden-1}) and we  obtain
\begin{equation*}
 \int_{B_R}\varphi_R \vert \vec{\nabla} \theta\vert^{2} dx=   \int_{B_R} \Delta \varphi_R  \left( \frac{\theta^2}{2}  \right)dx + \int_{B_R} \vec{\nabla} \varphi_R \cdot \left( \frac{ \theta^2}{2} \vu\right)dx. 
\end{equation*}
Moreover, since $\varphi_R (x)=1$ in the ball $B_{R/2}$,  $\text{supp}(\vec{\nabla} \varphi_R) \subset C(R/2,R)$ and $\text{supp}(\Delta \varphi_R) \subset C(R/2,R)$;   we have 
\begin{equation}\label{Estim-1}
\int_{B_{R/2}} \vert \vec{\nabla} \theta\vert^{2} dx \leq   \int_{C(R/2,R)} \Delta \varphi_R  \left( \frac{\theta^2}{2}  \right)dx + \int_{C(R/2,R)} \vec{\nabla} \varphi_R \cdot \left( \frac{ \theta^2}{2} \vu\right)dx. 
\end{equation}
In this estimate, we shall study each term on the right-hand side. For the first term, we use the H\"older inequalities with $\frac{1}{q}+\frac{2}{p}=1$ (hence we have $\frac{3}{q}-2=1-\frac{6}{p}$) and the fact that $\varphi_R(x)=\varphi(x/R)$ to write
\begin{equation}\label{Estim-2}
\begin{split}
\int_{C(R/2,R)} \Delta \varphi_R  \left( \frac{\theta^2}{2}  \right)dx \leq &\, C\, \left( \int_{C(R/2, R)} | \Delta \varphi_R  |^q dx  \right)^{1/q}\, \left( \int_{C(R/2, R)} | \theta |^p dx \right)^{2/p}\\
\leq & C\, R^{\frac{3}{q}-2} \| \Delta \varphi \|_{L^q(C(1/2,1))} \left( \int_{C(R/2, R)} | \theta |^p dx \right)^{2/p}\\
\leq & C\, R^{1-\frac{6}{p}}  \left( \int_{C(R/2, R)} | \theta |^p dx \right)^{2/p}. 
\end{split}
\end{equation}
For the second term, we use again the H\"older inequalities with $\frac{1}{r}+\frac{3}{p}=1$ (hence we have $\frac{3}{r}-1=2-\frac{9}{p}$) and with $\frac{3}{p}=\frac{2}{p}+\frac{1}{p}$  to write
\begin{equation}\label{Estim-3}
\begin{split}
 \int_{C(R/2,R)} \vec{\nabla} \varphi_R \cdot \left( \frac{ \theta^2}{2} \vu\right)dx \leq &\, C\, \left( \int_{C(R/2, R)} | \vec{\nabla} \varphi_R |^r dx  \right)^{1/r}\, \left( \int_{C(R/2, R)} | \theta^2 \, \vu  |^{p/3} dx  \right)^{3/p} \\
 \leq & \, C \, R^{\frac{3}{r}-1} \| \vec{\nabla} \varphi \|_{L^r(C(1/2,1))}\, \left( \int_{C(R/2,R)} | \theta |^p dx \right)^{2/p}\, \left( \int_{C(R/2,R)} | \vu  |^p dx \right)^{1/p} \\
 \leq & \, C \, R^{2-\frac{9}{p}}\, \left( \int_{C(R/2,R)} | \theta |^p dx \right)^{2/p}\, \left( \int_{C(R/2,R)} | \vu  |^p dx \right)^{1/p}. 
\end{split}
\end{equation}
Thus, the wished estimate (\ref{Caccioppoli-theta}) follows from the estimates (\ref{Estim-1}), (\ref{Estim-2}) and (\ref{Estim-3}). Proposition \ref{Caccioppoli} is proven. \finpv

\medskip 

With the Proposition \ref{Caccioppoli} at our disposal, we are able to proof the Theorem \ref{Th3}.  The main idea is to use this Caccioppoli-type estimate,  together with the information below on $\vu$ and $\theta$, to prove first that $\theta =0$. This identity then yields $\vu=0$ and $P=0$. 

\subsection{The case $\vu \in L^{\frac{9}{2}, q}(\Rt)$, $\theta \in  L^{\frac{9}{2},q}(\Rt)$  with $1\leq q<+\infty$  and $\vec{\bf g} \in \mathcal{W}^{1,9/2} \cap W^{1,\infty}(\Rt)$.}\label{Sec-Liouville}

With this information on $\vu$,  $\theta$ and $\vec{\bf g}$,  by Theorem \ref{Th2} we have $(\vu, P, \theta) \in \mathcal{C}^{2}(\Rt) \times \mathcal{C}^{1}(\Rt) \times \mathcal{C}^{2}(\Rt)$. Then, we set $p=\frac{9}{2}$ and by Proposition \ref{Caccioppoli} for all $R\geq 1$ we have the estimate 
\begin{equation}
\int_{B_{R/2}} | \vec{\nabla} \theta |^2 dx \leq  c\,   \left( \int_{C(R/2, R)} | \theta |^{\frac{9}{2}} dx \right)^{\frac{4}{9}} \left[  R^{-\frac{1}{3}}  +  \left( \int_{C(R/2,R)} | \vu  |^{\frac{9}{2}} dx \right)^{\frac{2}{9}}\right].
\end{equation} 
At this point, we recall the following estimate (see \cite[Proposition $1.110$]{DCh} for a proof): 
\begin{equation}
\int_{B_R} | f |^{\frac{9}{2}} dx \leq c\, R^{3(1-9/2q)} \| f \|^{\frac{9}{2}}_{L^{9/2, \infty}} \leq c\, R^{3(1-9/2q)} \| f \|^{\frac{9}{2}}_{L^{9/2, q}}.  
\end{equation}
We thus set $f= \mathds{1}_{C(R/2,R)} \theta$ and $f= \mathds{1}_{C(R/2,R)} \vu$; and we get back to the previous estimate to write
\begin{equation}
\begin{split}
\int_{B_{R/2}} | \vec{\nabla} \theta |^2 dx \leq &\, c \| \mathds{1}_{C(R/2,R)} \theta \|^{2}_{L^{9/2,q}} \left[ 1 + \| \vu \|_{L^{9/2,q}} \right]. 
\end{split}
\end{equation}
Thereafter, since $\theta \in L^{9/2,q}(\Rt)$ with $1 \leq  q <+\infty$ we can apply the dominated convergence theorem in the Lorentz spaces (see \cite[Theorem $1.2.8$]{DCh})  to obtain that $\ds{\lim_{R\to +\infty}\| \mathds{1}_{C(R/2,R)} \theta \|^{2}_{L^{9/2,q}} =0 }$. Consequently, in the last estimate we take the limit when $R\to +\infty$ to conclude that $\ds{\| \vec{\nabla} \theta \|_{L^2}=0}$. From this identity and the fact that $\theta \in L^{9/2,q}(\Rt)$ we finally get the identity $\theta =0$. 

\medskip

Once we have $\theta=0$, the coupled system (\ref{Boussinesq-Parabolic}) reduces to the stationary Navier-Stokes equation
\begin{equation}
-\Delta \vu + \text{div}(\vu \otimes \vu) + \vec{\nabla} P=0, \quad \text{div}(\vu)=0.
\end{equation}
Since $\vu \in L^{9/2,q}(\Rt)$ (with $9/2<q<+\infty$) is a smooth solution by \cite[Theorem $1$]{Jarrin} we obtain the identities $\vu=0$ and $P=0$.

\subsection{The case $\vu \in L^{p,\infty}(\Rt)$, $\theta \in L^{p,\infty}(\Rt)$ and $\vec{\bf g} \in \mathcal{W}^{1,p} \cap W^{1,\infty}(\Rt)$ with $3<p < \frac{9}{2}$.}

With this information on $\vu$ and $\theta$; and by Proposition \ref{Prop2} we have $\vu, \theta \in L^{\infty}(\Rt)$. Then, by Lemma \ref{Lem-Interpolation} and by the well-known embedding of the Lebesgue spaces into the Lorentz spaces we  obtain $\vu, \theta \in L^{9/2,q}(\Rt)$ with $9/2 < q < +\infty$. Thereafter, by the point \ref{Sec-Liouville} above we con directly conclude the identities $(\vu, P, \theta)=(0,0,0)$. 

\medskip 

Theorem \ref{Th3} is proven. \finpv

\section*{Author's declarations}
Data sharing not applicable to this article as no datasets were generated or analyzed during the current study. This work has not received any financial support. In addition, the author declares that he has no conflicts of interest.

\end{document}